%% file: FTHHH.tex
\documentclass[11pt]{amsart}
\pdfoutput=1
\usepackage{etex}
\usepackage[hmargin=3cm,vmargin=3cm]{geometry}

\usepackage[latin1]{inputenc}
\usepackage{xspace,amssymb,amsfonts,euscript}
\usepackage{amsthm,amsmath,amscd,stmaryrd,latexsym,youngtab,url}
\usepackage[nohug]{diagrams}
\usepackage{palatino, euscript}
\usepackage{ifpdf}

\usepackage[backend=bibtex,doi=false,isbn=false,url=false,style=alphabetic,maxnames=6]{biblatex}
\bibliography{mastercopy.bib}

\usepackage{graphicx}
\usepackage[all]{xy}
\SelectTips{cm}{}

\usepackage{tikz, tikz-cd}
\usetikzlibrary{matrix, arrows, calc}

\usepackage[dvips]{epsfig}
\usepackage{pinlabel}

\newcommand{\ig}[2]{\vcenter{\xy (0,0)*{\includegraphics[scale=#1]{diagrams/#2}} \endxy}}

\IfFileExists{comments.sty}{\usepackage{comments}}

\RequirePackage{color}
\definecolor{myred}{rgb}{0.75,0,0}
\definecolor{mygreen}{rgb}{0,0.5,0}
\definecolor{myblue}{rgb}{0,0,0.65}

  \RequirePackage[pdftex,
   colorlinks = true,
   urlcolor = myblue, 
   citecolor = mygreen, 
   linkcolor = myred, 
   ]
 {hyperref}

\newtheorem{thm}{Theorem}[section]
\newtheorem{lemma}[thm]{Lemma}
\newtheorem{theorem}[thm]{Theorem}

\newtheorem{prop}[thm]{Proposition}
\newtheorem{proposition}[thm]{Proposition}
\newtheorem{cor}[thm]{Corollary}

\newtheorem{conj}[thm]{Conjecture}
\newtheorem{conjecture}[thm]{Conjecture}

\newtheorem*{prop*}{Proposition}

\theoremstyle{definition}
\newtheorem{defn}[thm]{Definition}
\newtheorem{definition}[thm]{Definition}

\newtheorem{ex}[thm]{Example}
\newtheorem{example}[thm]{Example}

\theoremstyle{remark}
\newtheorem{remark}[thm]{Remark}
\newtheorem{observation}[thm]{Observation}

\numberwithin{equation}{section}

\newcommand{\ov}[1]{\overline{#1}}
\newcommand{\Sign}{\operatorname{sgn}}
\newcommand{\ga}{\gamma}
\newcommand{\w}{\omega}
\renewcommand{\t}{\tau}
\newcommand{\RHom}{\operatorname{RHom}}

    \def\QM{{\mathbb{Q}}}
    
    \def\SM{{\mathbb{S}}}

    \def\ZM{{\mathbb{Z}}}

    \def\AC{{\mathcal{A}}}
    
    \def\CC{{\mathcal{C}}}
    \def\DC{{\mathcal{D}}}

\def\HB{{\mathbf H}}

    \def\KC{{\mathcal{K}}}

    \def\PC{{\mathcal{P}}}

  \def\vb{{\mathbf v}}  
  \def\wb{{\mathbf w}}

\def\a{\alpha}
\def\b{\beta}
\newcommand{\be}{\beta}

\def\d{\delta}

\def\l{\lambda}

\def\s{\sigma}

\def\t{\tau}

\let\phi=\varphi
\let\tilde=\widetilde

\usepackage{bbm}
\def\C{{\mathbbm C}}

\def\Z{{\mathbbm Z}}
\def\Q{{\mathbbm Q}}
\def\1{\mathbbm{1}}
\newcommand{\one}{\1}

\newcommand{\un}{\underline}

\newcommand{\ot}{\otimes}

\newcommand{\co}{\colon}
\newcommand{\ip}[1]{\langle #1\rangle}

\renewcommand{\to}{\rightarrow}

\renewcommand{\sl}{\mathfrak{sl}}

\newcommand{\refequal}[1]{\xy {\ar@{=}^{#1}
(-1,0)*{};(1,0)*{}};
\endxy}

\newcommand{\Hom}{{\rm Hom}}

\newcommand{\END}{{\rm END}}

\newcommand{\op}{{\rm op}}

\newcommand{\Id}{\operatorname{Id}}

\newcommand{\Tr}{\operatorname{Tr}}

\newcommand{\Sh}{\operatorname{Sh}}
\newcommand{\inv}{^{-1}}
\newcommand{\vect}{\textbf{-vect}}
\newcommand{\In}{\operatorname{In}}
\newcommand{\Out}{\operatorname{Out}}
\newcommand{\Tw}{\operatorname{Tw}}

\newcommand{\Cone}{\textrm{Cone}}

\newcommand{\SBim}{\SM\textrm{Bim}}

\newcommand{\HHH}{\operatorname{HHH}}
\newcommand{\HH}{\operatorname{HH}}
\newcommand{\FHilb}{\operatorname{FHilb}}
\newcommand{\K}{\KC}

\newcommand{\FT}{\operatorname{FT}}
\newcommand{\HT}{\operatorname{HT}}
\newcommand{\Br}{\operatorname{Br}}

\begin{document}


\begin{abstract}
We introduce a new method for computing triply graded link homology, which is particularly well-adapted to torus links.  Our main application is to the $(n,n)$-torus links, for which we give an exact answer for all $n$.  In several cases, our computations verify conjectures of Gorsky \emph{et al} relating homology of torus links with Hilbert schemes.
\end{abstract}

\title{On the computation of torus link homology}

\author{Ben Elias} \address{University of Oregon, Eugene}

\author{Matthew Hogancamp} \address{Indiana University, Bloomington}

\maketitle

\tableofcontents

\input FTHHHintro.tex

\input FTHHHbackground.tex

\input FTHHHconvolve1.tex

\input FTHHHconvolve3.tex

\input FTHHHnumerology.tex

\input FTHHHcomputations.tex

\printbibliography

\end{document}

%% file: FTHHHintro.tex
\section{Introduction}
\label{sec-intro}

Triply-graded Khovanov-Rozansky homology is a link homology theory which was originally introduced by Khovanov-Rozansky \cite{KhoRoz08} using matrix factorizations, but was soon after
reinterpreted by Khovanov \cite{Khov07} using the Hochschild homology of Soergel bimodules. It has generated a great deal of interest, admitting spectral sequences which converge to
various $\sl_n$-link homology theories \cite{RasDifferentials}, and having deep connections to the representation theory of Hecke algebras in type $A$.

Khovanov's construction begins with a braid $\be$ on $n$ strands. To such a braid, Rouquier \cite{RouqBraid-pp} has associated a complex (up to homotopy) $F(\be)$ of Soergel bimodules,
which are certain graded bimodules \cite{Soer07} over the polynomial ring $R = R_n = \QM[x_1, \ldots, x_n]$ in $n$ variables. More precisely, Rouquier associates a complex to each braid
generator (e.g. over- or under-crossing). From this, one obtains a complex for any braid diagram by taking the tensor product of these elementary complexes. Rouquier proves that two
braid diagrams for the same braid yield complexes which are canonically isomorphic in the homotopy category of $R$-bimodules.

Khovanov \cite{Kh07} observed that taking the closure $\ov{\be}$ of a braid $\be$ should correspond to identifying the right and left actions of $R$, or rather the higher derived functors of this operation.  These higher derived functors are known as Hochschild homology and are denoted by $\HH_i(R;M)$; when $R$ is understood we write $\HH_i(M)=\HH_i(R;M)$.  Khovanov proved that the complex obtained by applying $\HH_i$ to each Soergel bimodule in a Rouquier complex $F(\be)$ yields a complex of vector spaces which (up to homotopy) depends only on the closure $\ov{\be}$, and thus the homology groups of this complex are link invariants of $\ov{\be}$.  The three gradings come from the Hochschild homological grading, the usual homological grading, and the internal grading of the Soergel bimodules.

It is well known (see \cite{Kh07} and references therein) that if $R=R_n$ is a polynomial ring and $M$ is an $R$-bimodule, then there is an isomorphism between Hochschild homology group $\HH_i(R;M)$, and the Hochschild \emph{cohomology} group $\HH^{n-i}(R;M)$.  The Hochshild cohomology groups are the the higher derived functors of $M\mapsto \Hom_{(R,R)}(R,M)$, the space of $R$-bimodule maps from the monoidal identity $R$.  For the remainder of this paper, we work exclusively with Hochschild cohomology.   

The Khovanov-Rozansky homology of torus links has deep connections to Hilbert schemes, rational Cherednik algebras, and refined Chern-Simons theory \cite{GORS,GorNeg15,GorNegRa-un}.    At the moment these connections are purely conjectural, but they suggest that the Khovanov-Rozansky homologies of torus links are quite interesting objects.  Up until now, however, the connection with other subjects has been difficult to verify, since the computation of Khovanov-Rozansky homology is quite challenging in practice.   In this paper we introduce a new method for computing Khovanov-Rozansky homology which seems particularly well adapted to compute homologies of torus links.  In particular, we provide a remarkably simple description of the triply-graded homology of the $(n,n)$ torus links, in Theorem \ref{introthm-FTpoincare}.

In \S \ref{subsec:magicFormula} we compare our results with the predictions of Gorsky-Negut \cite{GorNeg15} (also Gorsky-Negut-Rasmussen \cite{GorNegRa-un}) coming from flag Hilbert schemes; in every case we have checked, they match identically.  Previous checks of the connection with Hilbert schemes have been limited to the cases $n=2,3$, but with our method we are able to verify the predictions of \emph{loc.~cit.} for $n\leq 4$.

On the other hand, it is difficult to compare our results on $(n,n)$ torus links with the conjectures of Gorsky-Oblomkov-Rasmussen-Shende \cite{GOR} since they focus on the case of the $(n,m)$-torus \emph{knots}, that is, when $n$ and $m$ are coprime, and much less is known about the link case.  Nonetheless, from P.~Etingof we learned that the ring of $k$-quasi-invariants for $S_n$ acting on $\Q[x_1,\ldots,x_n]$ (see \cite{EtingStrick-pp} for a survey) is a representation of the rational Cherednik algebra for $\sl_n$, and is the correct replacement for the simple module $L_{m/n}$ which appears in \cite{GORS} when $m=kn$.  Thus, the minimal Hochschild degree part of the Poincar\'e series of the $(n,nk)$ torus links is expected to equal the Hilbert series for a certain filtration on the ring of $k$-quasi-invariants.  However, it is not clear how to filter the ring of quasi-invariants in an appropriate way, so we we will not say more about this connection in this paper.

In Appendix \ref{sec:appendix} we include some additional computations.  We found that the Poincar\'e polynomial of $\HH^0$ of the $(n,n+1)$ torus knot is given by the $q,t$ Catalan number for $n=2,3,4$, which verifies a conjecture of Gorsky's \cite{Gorsky12} for these knots.

Our particular interest in the homology of the $(n,n)$ torus links stems from the fact that this triply graded vector space parametrizes maps from the identity Soergel bimodule $R$ to the Rouquier complex associated to the full twist braid $\FT_n$.  The computation above is used in forthcoming work of the authors, in which we decompose the Soergel category into its ``eigencategories'' for the action of $\FT_n$, thereby laying the groundwork for the study of the categorical representation theory of Hecke algebras.

\subsection{Motivation from categorical representation theory}
\label{subsec:catRepThry}

It was shown by Khovanov-Thomas \cite{KT} that Rouquier complexes give a faithful action of the braid group on the homotopy category of Soergel bimodules. For this reason, the collection of Rouquier complexes is often refered to as a categorification of the braid group, but this is somewhat misleading, as this particular action of the braid group is intricately tied to its Hecke quotient. Soergel \cite{Soer07} proved that Soergel bimodules over $R_n$ categorify the Hecke algebra $\HB = \HB_n$ of the symmetric group $S_n$. Note that $\HB$ is linear over the ring $\ZM[Q,Q^{-1}]$, where $Q$ is categorified by the grading shift of an $R$-bimodule\footnote{Works of Soergel and those who followed him often work over the ring $\ZM[v,v^{-1}]$ instead; to compare conventions, use the equality $Q = v^{-1}$.}; for this to work correctly, $R$ is graded so that $\deg x_i = 2$ for all $1 \leq i \leq n$. Taking the image of a Rouquier complex in this Grothendieck group yields the familiar quotient map from the braid group (or its group algebra over $\ZM[Q,Q^{-1}]$) to the Hecke algebra. However, the fact that Rouquier complexes only reflect the ``Hecke-type'' actions of the braid group is an advantage, not a limitation, as one can lift results from the representation theory of
Hecke algebras to study the homotopy category of Soergel bimodules. This paper can be understood and appreciated without a foray into categorical representation theory, but we provide
some brief motivation here.

The Hecke algebra admits a \emph{sign representation} $\Sign_n$, on which each of its standard generators (the images of the overcrossings) acts by $-Q^{-1}$. The projection from an
arbitrary Hecke algebra representation to its isotypic component for the sign representation is an idempotent often known as a (generalized) Jones-Wenzl projector, after the
corresponding idempotent in the Temperley-Lieb algebra \cite{Jones01, Wenzl}. This projection can not be defined in $\HB$ itself, requiring certain scalars to be inverted (like
$Q+Q^{-1}$, for example). It can be defined in the base change $\HB \ot_{\ZM[Q,Q^{-1}]} \ZM((Q))$.

In \cite{Hoga-pp}, the second author constructs an infinite complex of Soergel bimodules $P_n$ which categorifies this Jones-Wenzl projector.  In this paper we study a finite complex $K_n$ which
categorifies the ``renormalized'' Jones-Wenzl projector, a rescaling of the projector which is actually defined within $\HB$ before base change. The fact that the Jones-Wenzl idempotent
projects to the sign representation is categorified by the fact that the Rouquier complex for an overcrossing, acting by tensor product on $K_n$, will simply act by a homological and a
grading shift.

The inductive construction of $K_n$ itself also is motivated by the representation theory of the Hecke algebra. When $\Sign_n$ is induced from $\HB_n$ to $\HB_{n+1}$, it splits into two
irreducible representations, $\Sign_{n+1}$ and another representation $V$. This splitting is actually the eigenspace decomposition for the Young-Jucys-Murphy operator $y_{n+1}$, a
certain element of the braid group on $n+1$ strands which commutes with any braid on the first $n$ strands. If $\ga$ is the eigenvalue corresponding to $V$, then $y_{n+1} - \ga$ kills
$V$, and thus is equal to the projection to $\Sign_{n+1}$ up to scalar.  If $k_n\in \HB_n\subset \HB_{n+1}$ denotes the renormalized projection onto the sign representation, then by the previous discussion there is a linear relation
\[
k_{n+1} = k_ny_{n+1} - \gamma k_n.
\]
On the categorical level, this relation becomes an exact triangle.  More precisely, there is a grading shift $\Gamma$ and a chain map $\phi:\Gamma K_n\rightarrow K_nF(y_{n+1})$ such that $K_{n+1}:=\Cone(\phi)$ categorifies the renormalized projection to $\Sign_{n+1}$. Recall that $F(y_{n+1})$ indicates the Rouquier complex associated to $y_{n+1}$.  It turns out that $\gamma = Q^2$, and $\Gamma=Q^2$ is simply the functor which shifts internal degree up by 2.   This chain map is constructed in \cite{Hoga-pp}, and we recall the basics in \S \ref{subsec-Hoga-pp}. 


This is an example of \emph{categorical diagonalization}, a concept which is developed in forthcoming work \cite{EHDiag}. The chain map $\phi$ mentioned above is an \emph{eigenmap}; in our categorification of various concepts in linear algebra, the cones of eigenmaps are used to categorify the operators $(A - \l I)$ for an eigenvalue $\l$ of an operator $A$. This makes the computation of $\HH^0$ of a braid particularly significant, because it describes the space of maps from the (shifted) monoidal identity, which are potential eigenmaps. In fact, the main result of this paper will be used as a lemma in \cite{EHDiag} to prove that the full twist in the braid group has enough eigenmaps and therefore is categorically diagonalizable. We use this to construct categorical projections to arbitrary irreducible representations of the Hecke algebra, not just the sign representation.

In this paper, our focus is on computation: the existence of $K_n$ is known by other means, and we use the recursive definition of $K_n$ to compute link invariants.  This strategy is outlined below.

\subsection{Our method, decategorified}
\label{subsec:methodsDecat}

The Hecke algebra $\HB_n$ is isomorphic to a quotient of the group algebra $\Z[Q,Q\inv][\Br_n]$ where we identify
\[
\ig{1}{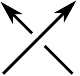}\ -\ \ig{1}{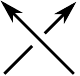}=(Q-Q\inv)\ig{1}{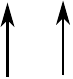}.
\]
The Jones-Ocneanu trace $\Tr:\HB_n\rightarrow \Z[Q^{\pm},A^{\pm}]$ is such that $\Tr(\b)$ is the Homfly polynomial of the braid closure $\hat{\b}$.  Using the skein relation above, and the formula which defines $\Tr$, one can in principal compute the Homfly polynomial for any link.  In \eqref{introeqn-Kn-cycles} we introduce another skein-like relation which is often useful.

There are elements $k_n\in \HB_n$ defined inductively by $k_1=1\in \HB_n$, and
\begin{equation} \label{introEqn-YJM-k}
\begin{minipage}{1in}
\labellist
\small
\pinlabel $k_{n}$ at 28 29
\endlabellist
\includegraphics[scale=1]{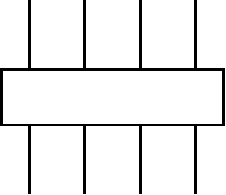}
\end{minipage}
\ \ = \ \ 
\begin{minipage}{1in}
\labellist
\small
\pinlabel $k_{n-1}$ at 27 39
\endlabellist
\includegraphics[scale=1]{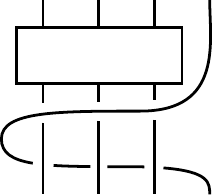}
\end{minipage}
\ \ - \ \ 
Q^2 \ \begin{minipage}{1in}
\labellist
\small
\pinlabel $k_{n-1}$ at 28 29
\endlabellist
\includegraphics[scale=1]{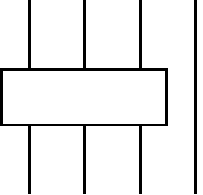}
\end{minipage}.
\end{equation}
The element $k_n$ is a renormalized projection onto the sign representation: $k_n \sigma_i^\pm = -Q^{\mp} k_n = \sigma_i^\pm k_n$ for all $1\leq i\leq n-1$, where $\sigma_i$ denotes the elementary braid generator.  This implies that
\begin{equation} \label{introeqn-Kn-cycles}
\begin{minipage}{1in}
\labellist
\small
\pinlabel $k_{n-1}$ at 28 32
\endlabellist
\includegraphics[scale=1]{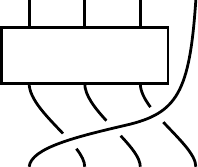}
\end{minipage}
\ \ = \ \ 
(-Q)^{n-1}\begin{minipage}{1.08in}
\labellist
\small
\pinlabel $k_{n}$ at 27 24
\endlabellist
\includegraphics[scale=1]{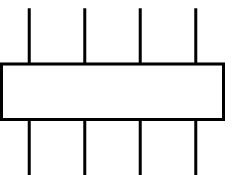}
\end{minipage}
 \ \ + \ \ 
Q^2\:\begin{minipage}{1.1in}
\labellist
\small
\pinlabel $k_{n-1}$ at 28 32
\endlabellist
\includegraphics[scale=1]{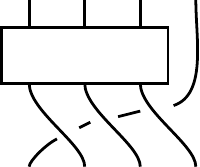}
\end{minipage}.
\end{equation}
The Jones-Ocneanu trace $\Tr(k_n)$ can easily be computed inductively from this equation using the invariance of $\Tr$ under the Markov moves together with the fact that $k_n$ absorbs crossings.  We can use $k_n$ to compute link invariants:
\begin{itemize}
\item Assume $\b$ is given.  Choose a crossing $x$ in $\b$.  Place $k_1$ somewhere in the vicinity of this crossing.  Since $k_1$ is the idenity element of $\HB_1$ , this does not change the element $\b\in \HB_n$.
\item Apply the relation \eqref{introeqn-Kn-cycles}; one of the terms will involove the switched crossing $x\inv$, and in the other term $k_1$ will have grown to $k_2$, which now has the potential to absorb some crossings.  
\item Repeat.  That is, assume that $\b$ is a braid with a $k_\ell$ inserted somewhere.  After manipulating the diagram, if necessary, arrange the picture so that Equation \eqref{introeqn-Kn-cycles} can be applied.  In one of the resulting terms, some crossings will be switched, which in good situations will simplify $\b$.  In the other term, $k_\ell$ grows in size and can now absorb more crossings, also resulting in a simpler diagram.
\end{itemize}
If one is lucky, this process can be repeated until the trace $\Tr$ of the resulting terms is trivial to compute.  Torus links seem especially well adapted to the application of this trick. 

\begin{example}
Let $x=\sigma_1\in \HB_2$ denote the crossing.  The trefoil is the $(2,3)$ torus knot, and can be presented as the closure of $x^3$.   Equation \eqref{introeqn-Kn-cycles} says that $x=-Q k_2 + Q^2 x\inv$.  Multiplying by $x$ gives
\[
x^2 = k_2 + Q^2.
\]
Multiplying by $x$ again gives
\[
x^3 = (-Q\inv)k_2 + Q^2 x.
\]
The trace of $k_2$ is easy to compute, and the trace of $x$ is the Homfly polynomial of the unknot, up to normalization.  Thus, the trace of $x^3$ is expressed in terms of known quantities.  We can continue in this manner, obtaining
\[
x^{2m} = (Q^{2(1-m)} + Q^{2(2-m)}+\cdots + Q^{2(m-1)})k_2 + Q^{2m}
\]
and
\[
x^{2m+1} = (-Q\inv)(Q^{2(1-m)} + Q^{2(2-m)}+\cdots + Q^{2(m-1)})k_2 + Q^{2m}x,
\]
from which the Homfly polynomials of the $(2,m)$ torus links are readily computed.
\end{example}

\subsection{Our method, categorified}
\label{subsec:methodsCat}
In this paper we categorify the method outlined in the previous section.  As alluded to earlier in this introduction, the element $k_n\in \HB_n$ gets replaced by a finite complex $K_n$ of Soergel bimodules, and the relations \eqref{introEqn-YJM-k} and \eqref{introeqn-Kn-cycles} become exact triangles.  More precisely, there is a chain map constructed in \cite{Hoga-pp} from $Q^2K_{n-1}\rightarrow K_{n-1}F(y_{n+1})$, and $K_{n+1}$ is defined to be the mapping cone on this map.  The fact that $k_n$ absorbs crossings becomes the fact that $K_nF(\sigma_i)\simeq TQ\inv K_n \simeq F(\sigma_i)K_n$.  Here and throughout we use $T$ and $Q$ to denote the functors which increase homological degree and bimodule degree respectively.  We have, for instance, an equivalence
\begin{equation} \label{introeqn-Kn-cycles-Cat}
\begin{minipage}{1in}
\labellist
\small
\pinlabel $K_{n-1}$ at 28 32
\endlabellist
\includegraphics[scale=1]{diagrams/KrecursionRight}
\end{minipage}
\ \ \simeq \ \ \left(
(TQ\inv)^{1-n}\begin{minipage}{1.08in}
\labellist
\small
\pinlabel $K_{n}$ at 27 24
\endlabellist
\includegraphics[scale=1]{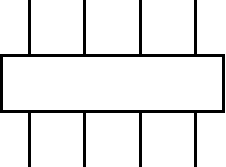}
\end{minipage}
 \ \ \longrightarrow  \ \ 
Q^2\ \begin{minipage}{1.1in}
\labellist
\small
\pinlabel $K_{n-1}$ at 28 32
\endlabellist
\includegraphics[scale=1]{diagrams/KrecursionLeft}
\end{minipage}\right),
\end{equation}
where the notation $A\simeq (B\rightarrow C)$ means that there is an exact triangle
\[
C\rightarrow A \rightarrow B\rightarrow T\inv C.
\]
The distinguished triangle \eqref{introeqn-Kn-cycles-Cat}  will be essentially the only weapon we need to attack our computations.  Suppose $F(\b)$ is a Rouquier complex that we would like to study.  Iterated application of the above exact triangle results in a certain kind of filtered complex (a convolution of a twisted complex; see below) which is homotopy equivalent to $F(\b)$, and whose subquotients are tensor products of Rouquier complexes and some $K_{\ell}$.  In favorable situations these have Hochschild cohomologies which are easy to compute.  There then arises the problem of recovering the homology of the total complex from the homology of its constituents.  For certain computations, we will see that this very serious complication is nullified by an equally serious miracle: the miracle of parity.

We first explain what sorts of filtered complexes we will use.  Suppose $A_i$ ($i\in I$) is a family of complexes, indexed by a finite partially ordered set $I$.  Suppose $d_{ij}:A_j\rightarrow A_i$ are a collection of linear maps such that
\begin{itemize}
\item $d_{ij}$ increases homological degree by 1.
\item $d_{ii}$ is the given differential on $A_i$.
\item $d_{ij}=0$ unless $i\ge j$.
\item the \emph{total differential} $d_{\text{tot}}:=\sum_{i\geq j}d_{ij}$ satisfies $d_{\text{tot}}^2=0$.
\end{itemize}
Then $C:=(\bigoplus_{i\in I} A_i, d_{\text{tot}})$ is a chain complex, which we call a \emph{convolution} of the $A_i$.  More precisely, this is the convolution of a one-sided twisted complex; see \cite{BonKapEnhanced} for more details on this construction in homological algebra.  We will also say that $C=\bigoplus_i A_i$ with \emph{twisted differential}, to indicate that the differential is not merely the sum of the differentials on the $A_i$.  The differential on this complex may be quite complicated, possibly sending terms in homological degree $m$ inside $A_j$ to terms in homological degree $m+1$ inside many different $A_i$, but it does respect the  order on $I$.  Thus, convolutions can also be thought of as certain kinds of filtered complexes, whose subquotients are the $A_i$.

Note that any exact triangle
\[
A_2\rightarrow C\rightarrow A_1 \rightarrow T\inv A_2
\]
gives rise to an equivalence $C\simeq (A_1\oplus A_2)$ with twisted differential.  Iterated mapping cones can be regarded as convolutions in a similar way.

Our main application is to the Rouquier complex $\FT_n$ associated to the full twist braids. In \S \ref{sec-resolve1} we iterate the equivalence \eqref{introeqn-Kn-cycles-Cat}, obtaining a convolution description of $F(y_n)$, the Rouquier complex associated to the Young-Jucys-Murphy braid.   Using the relation $\FT_n = \FT_{n-1}F(y_n)$, we then prove:
\begin{theorem}\label{introThm:FTconv}
We have $\FT_n \simeq \bigoplus_v q^{k} D_v$ with twisted differential.  The sum is over sequences $v\in \{0,1\}^n$ such that $v_n=1$.  Here, $q = Q^2$ indicates a grading shift, $k$ is the number of zeroes in $v$, and $D_v$ is described below.  The differential respects the anti-lexicographic order on sequences.
\end{theorem}
In the antilexicographic order we regard $(\ast,1)$ as larger than $(\ast,0)$, where $\ast$ denotes any sequence of zeroes and ones.  For each sequence $v \in \{0, 1\}^n$---which we will call a shuffle---there is a complex which we call $D_v$.  For example, here is $D_{10101101}$, which occurs (up to shift) in the expression of $\FT_8$.
\[
D_{10101101}=
\begin{minipage}{1.1in}
\labellist
\small
\pinlabel $\FT_3$ at 18 40
\pinlabel $K_5$ at 65 40
\endlabellist
\begin{center}\includegraphics[scale=1]{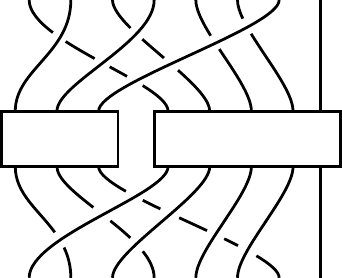}\end{center} 
\end{minipage}
\]
Inside $v$, the zeroes indicate which strands are connected to the full twist $\FT_k$, and the ones indicate which are connected to $K_\ell$, for $k + \ell = n$.   

Then, of course, one wants to compute the Hochschild homology of the complexes $D_v$.  Let us be precise.  Given a complex $F$ of Soergel bimodules, let $\HH^i(C)$ denote the complex obtained by applying the functor $\HH^i$ to each bimodule, and let $\HH(C) = \oplus \HH^i(C)$. Let $\HHH(C)$ denote the cohomology of the complex $\HH(C)$.  

Because Hochschild cohomology of a complex $C$ is unchanged by conjugation $C\mapsto FCF\inv$ for any invertible complex $F$, we can move part of $D_v$ from the bottom to the top, yielding the complex $C'_v$:
\[
C'_{10101101} \ = \ \begin{minipage}{1.5in}
\labellist
\small
\pinlabel $\FT_4$ at 20 20
\pinlabel $K_3$ at 70 20
\endlabellist
\begin{center}\includegraphics[scale=1]{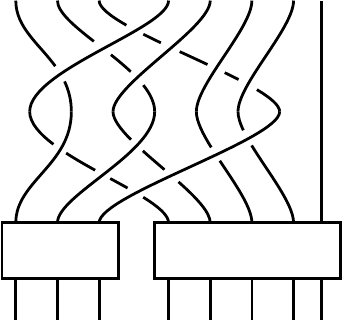}\end{center} 
\end{minipage}
\]
Note that $\FT_n = C'_{00\cdots 0}$.

For purely combinatorial reasons, we work instead with a similar complex $C_v$, which is defined by the same expression as $C_v'$, but with $K_\ell$ replaced by its reduced version $\hat{K}_\ell$.  Reduced complexes are discussed in \S\ref{subsec-reducedcomplexes}. The effect this has on Poincare polynomials is multiplication by a factor of $(1-Q^2)$.  

Let $v \cdot w$ denote the concatenation of two shuffles (sequences of zeroes and ones). For any shuffle $v$, we can use our distinguished triangle for $K_n$ to prove the following:

\begin{proposition}
We have $\HH(C_{v \cdot 0})\simeq \Big(\HH(C_{1 \cdot v})\rightarrow Q^2\HH(C_{0 \cdot v})\Big)$.
\end{proposition}
Next, we can use some relatively easy arguments involving the complex $K_n$ to prove that $\HH(C_{v \cdot 1})$ is just a direct sum of shifted copies of $\HH(C_v)$. For readers familiar with knot theory, this last statement should be thought of as analogous to the Markov move; it allows us to reduce the number of strands by $1$. Finally, a simple observation (pertaining to reduced complexes) allows one to replace the computation of $\HH(C_{000\cdots 0})$ with $\HH(C_{100\cdots 0})$. Combining these three operations, we obtain a recursive convolution description of any $\HH(C_v)$. This is the main result of \S\ref{subsec-newRecursion}.

Let us return to the computation of the Hochschild cohomology of the Rouquier complex for the full twist $\FT_n$ on $n$ strands.  


We are interested in the cohomology $\HHH(C_v)$ of the complexes $\HH(C_v)$. However, in general, the cohomology of a convolution of complexes is not the direct sum of the cohomology of the individual complexes; instead, there is a spectral sequence relating the two. Our final argument comes from observing a parity miracle! We prove inductively that $\HHH(C_v)$ is concentrated in even homological degrees. This forces every spectral sequence in sight to degenerate at the $E_1$ page, and implies that our convolution description of $\HH(C_v)$ gives a direct sum description of $\HHH(C_v)$. See Theorem \ref{thm:parityAndRecursion} and its proof for further discussion of this parity argument.

Thus, we have a recursive formula for the triply graded cohomologies $\HHH(C_v)$, and as a special case, a formula for $\HHH(\FT_n)$. We discuss this formula in the next section.


Let us pause to point out one of the subtleties we have ignored above. One can conjugate a complex by a braid and obtain a non-isomorphic complex with the same Hochschild cohomology. We
begin to apply this operation freely in \S\ref{sec-mainResult}. Above, we have stated that $\HH(C_{v \cdot 0})$ is a convolution of $\HH(C_{0 \cdot v})$ and $\HH(C_{1 \cdot v})$, but the
same statement does not hold for the original complexes. Instead, some conjugate of $C_{v \cdot 0}$ is a convolution of a conjugate of $C_{0 \cdot v}$ and a conjugate of $C_{1 \cdot v}$
(conjugating by different braids for each term). For purposes of Hochschild cohomology, this imprecision is harmless. However, were one to try to actually construct a chain map from $R$
to $C_v$ using this computation of $\HHH^0(C_v)$, then one would need to keep track of conjugation more carefully, which would be rather difficult.

On the other hand, the work done in \S\ref{sec-resolve1} describes full twist as a genuine convolution of complexes $D_v$, not complexes up to conjugation. This result is not actually
needed or used in the recursive computation of $\HHH(C_v)$ which is our main result. We include this auxiliary result because it can be used to construct chain maps from $R$ to $\FT_n$.
Our main theorem implies that the complexes $D_v$ satisfy a parity condition, and therefore our convolution description of $\FT_n$ induces a direct sum decompositions on Hochschild
cohomology. In particular, any chain map (up to homotopy) from $R$ to the complex $D_v$ (an element of $\HHH^0(D_v)$) can be extended uniquely (up to homotopy) to a chain map from $R$ to
the entire complex $\FT_n$. We use this fact to construct eigenmaps to the full twist in \cite{EHDiag}.

In addition, the convolution description involving $D_v$ from \S\ref{sec-resolve1} can be adapted to other torus links, whereas the results \S \ref{sec-mainResult} are fundamentally tied
to the case of $(n,n)$-torus links.

\begin{remark} Suppose that one were interested in computing the Hochschild cohomologies of $\FT_n^{-1}$. One can produce a convolution description of $\FT_n^{-1}$ similar to the description of $\FT_n$ above, but the parity
miracle no longer holds! The corresponding spectral sequence is far from degenerate, and consequently there are extremely few (non-nulhomotopic) chain maps from $R$ to $\FT_n^{-1}$. This lack of symmetry between $\FT_n$ and $\FT_n^{-1}$ is an interesting and complicating feature in categorical representation theory.  \end{remark}

\subsection{The recursive formula}
\label{subsec:introFT}
We will find it convenient to use a non-standard choice of variables for our Poincare series. We let $t=T^2Q^{-2}$, $q=Q^2$, and $a = Q^{-2} A$, where $T$ denotes the usual homological
degree, $Q$ the bimodule degree (also called internal degree, or quantum degree), and $A$ the Hochschild degree. For instance, the Poincar\'e series of the polynomial ring $R=\Q[x_1,\ldots,x_n]$ is written $1/(1-q)^n$, and the Poincar\'e series of its Hochschild cohomology is $(1-q)^{-n}(1+a)^n$. In general, all our Poincar\'e series will be power series in the variables $q$, $a$, and $t^{\frac{1}{2}}$.

\begin{prop}\label{prop:fv-intro}
There is a unique family of polynomials $f_v(q,a,t)$, indexed by integers $n \ge 0$ and binary sequences $v\in \{0,1\}^n$, satisfying $f_{\emptyset}=1$ together with
\begin{subequations}
\begin{equation} \label{eq:newRecursion1}f_{v\cdot 1}(q,a,t) = (t^{|v|}+a )f_v\end{equation}
\begin{equation}\label{eq:newRecursion2}f_{v\cdot 0}(q,a,t) = q f_{0\cdot v} + f_{1\cdot v}\end{equation}
\end{subequations}
\end{prop}

The following theorem, together with the fact that the $f_v$ are rational functions in $q,a,t$ rather than $q,a,t^{\frac{1}{2}}$ implies that the parity miracle holds:

\begin{theorem}\label{introthm-FTpoincare}
The Poincar\'e series of $\HHH(C_v)$ is $f_v(q,a,t)$, where $C_v$ are the complexes from Definition \ref{def:Chat}. In particular, the Hochschild cohomology $\HHH(\FT)$ is given by $f_v(q,a,t)$ for $v = (00\cdots 0)$.  These homologies are all supported in even homological degrees.
\end{theorem}

These results are restated and proved in \S\ref{subsec:mainResult}. The proof of Theorem \ref{introthm-FTpoincare} comes from the convolution description of $\HH(C_v)$ discussed in the
previous section.

For the reader's edification, here are the complete power series for $\HHH(\FT_n)$ for $n=1,2,3$.
\begin{subequations}
\begin{equation} \nonumber f_0(q,a,t) = \frac{1+a}{1-q}. \end{equation}
\begin{equation} \nonumber f_{00}(q,a,t) = \frac{1+a}{(1-q)^2}(q+t-qt+a). \end{equation}
\begin{eqnarray} \nonumber f_{000}(q,a,t) & = &  \frac{1+a}{(1-q)^3}\Big((t^3q^2 + q^3t^2 - 2t^2q^2 - 2tq^3 - 2 qt^3 + t^3 + q^3 + tq^2 + qt^2+tq) \\ \nonumber && +(t^2q^2 - 2tq^2 - 2qt^2 + t^2 + q^2 + tq + t + q)a + a^2\Big). \end{eqnarray}
\end{subequations}

The recursion  can be unraveled into the equivalent recursion below, which is more complicated but faster to implement.

\begin{definition} \label{defn-f_v}
For each integer $n \ge 0$, we let $[n]:=\{1,\ldots,n\}$.  We identify subsets $v\subset [n]$ with binary sequences $v\in\{0,1\}^n$.  For each such $v\subset [n]$, we define a rational function $f_v(q,a,t)$ by $f_\emptyset=1$, together with the following rules:
\begin{subequations}
\begin{equation}\label{eq:oldRecursion1}f_{000 \cdots 0}(q,a,t) = \frac{1}{1-q}f_{100 \cdots 0}(q,a,t).\end{equation}
\begin{equation}\label{eq:oldRecursion2}f_{11 \cdots 1}(q,a,t) = \prod_{i=1}^n (t^{i-1} +  a),  \ \ \ \ \ \ \text{($n$ indices)}.\end{equation}
\begin{equation}\label{eq:oldRecursion3}
f_v(q,a,t) = \sum_{w\subset [k]} \PC_{v,w}(a,t)  q^{k-|w|} f_w(q,a,t),  \ \ \ \ \ \ \text{$v\neq (0\cdots 0)$ and $v\neq (1\cdots 1)$}.
\end{equation}
\end{subequations}
Here $k=n-|v|$ is the number of zeroes in $v$.  Definition \ref{defn-gradings} contains the description of $\PC_{v,w}(a,t)$, which is a product of $n-k$ factors each of the form $(t^{\ell+m}+ a)$ for various numbers $\ell$ and $m$ depending on the sequences $v$ and $w$.
\end{definition}

The equivalence between these recursions is proven in chapter \ref{sec-numerology}, which contains various such numerological considerations. For example, one can show that both rule \eqref{eq:oldRecursion1} and \eqref{eq:oldRecursion2} are actually just consequences of rule \eqref{eq:oldRecursion3} when applied verbatim, although this is not obvious. 

\begin{remark} In our original version of this manuscript, our convolution argument categorified the recursive formula of Definition \ref{defn-f_v} rather than Proposition
\ref{prop:fv-intro}, thus proving that this complicated recursion does compute $\HHH(C_v)$. Then we discovered the simpler recursion of Proposition \ref{prop:fv-intro}, drastically
simplifying our arguments. 
\end{remark}

Note that the contribution to higher Hochschild gradings comes only from a factor in rule \eqref{eq:oldRecursion2} and the factor $\PC_{v,w}(a,t)$ in rule \eqref{eq:oldRecursion3}; both of these become explicit monomials in $t$ upon setting $a=0$.  Thus, to understand the zero-th Hochschild degree coefficient, i.e. the polynomial $f_v(q,0,t)$, one use a simplified recursion relation.  Using this, we prove the following closed formula for the power series $f_{00\cdots 0}(q,0,t)$, also known as the Poincar\'e series of the zero-th Hochschild cohomology $\HHH^0(\FT)$:

\begin{theorem} \label{thm-degreezerosolve-intro}
The Hochschild degree zero part of the unreduced triply graded homology of $(n,n)$ torus links has Poincare series equal to
\[
F_n(q,t)=\sum_{\sigma}t^{a(\sigma)+b(\sigma)}q^{c(\sigma)}
\]
where the sum is over functions $\sigma:\{1,\ldots,n\}\rightarrow \Z_{\geq 0}$, and the integers $a(\sigma)$, $b(\sigma)$, $c(\sigma)$ are defined by
\begin{enumerate}
\item $a(\sigma)=\sum_{k\geq 0}\binom{|\sigma\inv(k)|}{2}$
\item $b(\sigma)$ is the number of pairs $(i,j)\in \{1,\ldots,n\}$ such that $i<j$ and $\sigma(j)=\sigma(i)+1$.
\item $c(\sigma)=\sum_{i=1}^n\sigma(i)$.
\end{enumerate}
\end{theorem}

\begin{example}
In case $n=1$ we have $F_1(q,t)=1+q+q^2+\cdots = 1/(1-q)$.
\end{example}

\begin{example}
In case $n=2$, $F_2(q,t)$ is the sum of monomials appearing in the following diagram:
\begin{diagram}[small]
t && tq && q^2 && q^3 && \cdots \\
q && tq^2 && tq^3 && q^4 && \cdots \\
q^2 && q^3 && tq^4 && tq^5 && \cdots \\
q^3 && q^4 && q^5 && tq^6 && \cdots \\
\vdots && \vdots && \vdots && \vdots &\ddots& 
\end{diagram}
After rearranging, this becomes $F_2(q,t)=t/(1-q) + q/(1-q)^2$, agreeing with $f_{00}(q,0,t)$ which was computed above.
\end{example}

The proof of this closed formula from Theorem \ref{introthm-FTpoincare} is a simple combinatorial argument, and is found in \S\ref{sec-numerology}. Unfortunately, the polynomials $\PC_{v,w}(q,a,t)$ are sufficiently complicated so that we have been unable to produce a closed formula for the higher Hochschild degrees along these lines.

We conclude with a recursion for the normalized polynomials $\tilde{f}_v(q,a,t):=(1-q)^{k}f_v(q,a,t)$, where $k$ is the number of zeroes in $v$.  The recursion of Proposition \ref{prop:fv-intro} immediately gives rise to
\begin{subequations}
\begin{equation} \label{eq:newNormRecursion1}\tilde{f}_{v\cdot 1}(q,a,t) = (t^{|v|}+a )\tilde{f}_v\end{equation}
\begin{equation}\label{eq:newNormRecursion2}\tilde{f}_{v\cdot 0}(q,a,t) = q \tilde{f}_{0\cdot v} +(1-q) \tilde{f}_{1\cdot v}\end{equation}
\end{subequations}

\begin{remark}\label{rmk:dealer} A clumsy card dealer has a deck of $n$ cards, some face up and some face down. When the dealer encounters a face down card, he deals it.  When the dealer encounters a face
up card, he puts it back on the bottom of the deck, sometimes remembering to flip it face down. Eventually, the deck is dealt (with probability $1$).  Every time the dealer deals a card, you, the player, choose whether to receive 1 silver coin, or a number of dollars equal to the number of face-up cards in the deck.   Then the coefficient of $a^kt^\ell$ in $\tilde{f}_v(q,a,t)$ is the number of ways of ending up with $k$ silver coins and $\ell$ dollars, weighted by their probability of occuring.  In particular, the coeffcient on $a^n$ is 1.  \end{remark}

Computer experiments suggest the following:

\begin{conj}\label{conj:symmetry}
We have the following symmetry: $\tilde{f}_{00\cdots 0}(q,a,t)=\tilde{f}_{00\cdots 0}(t,a,q)$.
\end{conj}

This symmetry would follow from a formula of Gorsky-Negut-Rasmussen, which we discuss now.

\subsection{Flag Hilbert schemes and a magic formula}
\label{subsec:magicFormula}
According to the remarkable work of Gorsky, Negut, and Rasmussen \cite{GorNeg15,GorNegRa-un}, triply graded link homology can be extracted from flag Hilbert schemes.\footnote{We warn the reader that this story is related to, but quite different from, other connections between link homology and algebraic geometry.}   Roughly, the picture looks like this: there is a space $\FHilb_n(\C^2)$, which parametrizes flags of ideals $I_1\subset \cdots \subset I_n\in \C[x,y]$ such that $I_{i}/I_{i-1}$ is 1-dimensional.   Associated to each $n$-strand braid, Gorsky-Negut-Rasmussen conjecture that there exists a line bundle (or sheaf, or complex of sheaves) on $\FHilb_n(\C^2)$ whose space of global sections recovers $\HHH^0(F(\b))$. The sheaves are meant to be equivariant with respect to an obvious action of $\C^\ast\times \C^{\ast}$, and the variables $q,t$ correspond to weights with respect to this action.  The variable $a$ can also be accounted for with more work.  A combination of Atiyah-Bott localization and careful analysis of the flag Hilbert scheme near its torus fixed points yields a remarkably simple combinatorial formula which, given their conjecture, will describe the knot homology of positive torus links.  Now we state the formulas, and we will make no further mention of the geometric foundations which motivate them.

Let $\lambda$ be a Young diagram, drawn in the ``English style'' as in:
\[
\ig{.6}{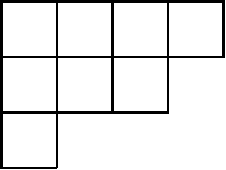}
\]
Suppose a box $c$ is in the $i$-th column and $j$-th row.  Here columns and rows are counted left-to-right and top-to-bottom, starting at zero.  To such a box we associated the monomial $z_c = t^i q^j$.\footnote{ In other variables, $z_c = T^{2i} Q^{-2x(c)}$, where $x(c)$ is the \emph{content} of the box $c$, which is $i-j$.}  To a Young diagram, we let $z_\lambda = z_\lambda(q,t)$ be the product of $z_c$ as $c$ ranges over all the boxes of $\lambda$. Note that $z_\lambda(q,t) = z_{\lambda^t}(t,q)$, where $\lambda^t$ is the transposed partition.

A box in $\lambda$ is \emph{removable} if $\lambda \smallsetminus c$ is a Young diagram. Now we discuss boxes, i.e. coordinates $(i,j)$, which need not be in the given Young diagram $\l$. We call $c \notin \lambda$ an \emph{outer corner} of $\lambda$ if the top left corner of $c$ coincides with the bottom right corner of a removable box in $\lambda$. We call a box $c \notin \lambda$ an \emph{inner corner} of $\lambda$ if $\lambda \cup c$ is a Young diagram. Pictorially, we have
\[
\ig{.6}{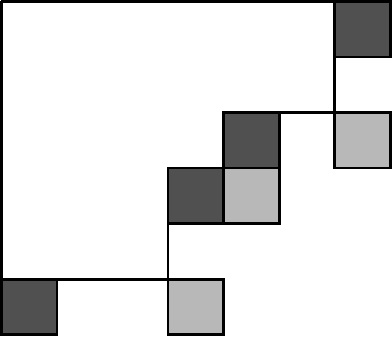}
\]
The inner corners are darkly shaded, and the outer corners are lightly shaded.  Let $\In(\lambda)$ and $\Out(\lambda)$ denote the sets of inner and outer corners of $\lambda$.  If $c\in \In(\lambda)$, then we define
\[
f_{\lambda,c}(q,t): = \frac{\prod_{d\in \Out(\lambda)}(z_c-z_d)}{\prod_{e\in \In(\lambda)\smallsetminus \{c\}}(z_c-z_e)}.
\]
It is easy to observe that $f_{\lambda,c}(q,t) = f_{\lambda^t, c^t}(t,q)$, where $c^t$ is the corresponding transposed inner corner of $\lambda^t$.

\begin{remark} Let $\{x_1, \ldots, x_n\}$ and $\{y_1, \ldots, y_{n+1}\}$ be two families of abstract variables. It is not a hard excercise to show that
\[
\sum_{i=1}^{n+1} \frac{\prod_j (y_i - x_j)}{\prod_{k \ne i} (y_i - y_k)} = 1.
\]
Applying this general formula to the definition above, one obtains
\begin{equation} \label{sumis1} \sum_{c \in \In(\lambda)} f_{\lambda, c} = 1, \end{equation}
a fact which has nothing to do with the combinatorics of partitions.
\end{remark}

A standard tableau $T$ can be thought of as a sequence of Young diagrams $T=(\lambda_1,\ldots,\lambda_n)$ such that $\lambda_1 = \square$ and $\lambda_{i+1}\smallsetminus \lambda_i = \square$. Let $\Sh(T) = \l_n$ denote the shape of $T$. To each tableau we set
\[
f_T(q,t):=\prod_{i=1}^{n-1} f_{\lambda_i,c_i}(q,t)
\]
where $c_i\in \In(\lambda_i)$ is the box such that $\lambda_{i+1}=\lambda_i\cup c_i$. Once more, $f_T(q,t) = f_{T^t}(t,q)$, where $T^t$ is the transposed tableau.

Finally, associated to a Young diagram $\lambda$, let $g_\lambda(q,a,t)$ denote the product over all boxes $c\in \lambda$ of $(1+az_c\inv)$.

\begin{conjecture}[Magic formula]\label{conj-magicFormula}
Let $F_{n,r}(q,a,t)$ denote the Poincar\'e series of $\HHH^0(\FT_n^{\otimes r})$.  Then
\begin{equation}\label{eq:magicFormula}
(1-q)^n F_{n,r}(q,a,t) = \sum_{T} z_{\Sh(T)}^r g_{\Sh(T)} f_T,
\end{equation}
a sum over all tableaux with $n$ boxes. In particular, the right hand side is symmetric under replacing $q$ with $t$, and thus so is the left hand side. We remind the reader that $z_{\Sh(T)}$ and $f_T$ are functions of $q$ and $t$, while $g_{\Sh(T)}$ is a function of $q$, $a$, and $t$.
\end{conjecture}

When $r=0$, the formula yields $1=\sum_T  g_{\Sh(\lambda)}(q,a,t)f_T(q,t)$, whose $a$-degree zero part follows from \eqref{sumis1}. We have verified the magic formula for $r=1$ and for $1\leq n\leq 4$, using a mathematica notebook which is available on the second author's website \cite{HogancampIUsite}. 

According to the magic formula, the $a$-degree $n$ part of the Poincar\'e series of the $(n,n)$ torus link is supposed to be $(1-q)^{-n}$ times
\[
\sum_T z_{\Sh(T)} z_{\Sh(T)}\inv f_T(q,t) = 1.
\]
The factor $z_{\Sh(T)}\inv$ comes from taking the $a$-degree $n$ part of $g_{\Sh(T)}(q,a,t)$.  This instance of the magic formula can be prove directly from our recursive description of this series (see Remark \ref{rmk:dealer}).  

Another consequence of the magic formula concerns the sub-maximal part of the Poincar\'e series of full twists:
\begin{conjecture}\label{conj:submaximal}
The $a$-degree $n-1$ part of the Poincar\'e series of $\HHH(\FT_n)$ is a geometric progression
\[
\frac{1}{(1-q)^n}\frac{1-(q+t-qt)^n}{(1-q)(1-t)} = \frac{1}{(1-q)^n}\Big(1+(q+t-qt)+\cdots + (q+t-qt)^{n-1}\Big)
\]
\end{conjecture}
We expect that this is not difficult to prove, but we do not do so here.  We have verified this conjecture up to $n=7$ using computer calculations.



In the algorithm to compute $f_{00\cdots 0}(q,a,t)$ using \eqref{eq:oldRecursion3}, one travels from the zero sequence $(00\cdots 0) \in \{0,1\}^n$ to the sequence $\emptyset \in
\{0,1\}^0$ by repeatedly choosing subsets (i.e. smaller sequences in $\{0,1\}^k$) of the previous set of zeroes. Our instinct indicates that such a sequence of sequences can be thought
of as encoding the entries in a Robinson-Shensted row-bumping algorithm, and can thus be assigned a tableau. The contribution to $f_{00\cdots 0}(q,a,t)$ coming from this sequence of
sequences and the contribution to $\sum_T z_{\Sh(T)} g_{\Sh(T)} f_T$ coming from the tableau $T$ have many superficial similarities, but no direct relation has yet been
found.

\subsection{Organization of the paper}

In \S\ref{sec-background} we provide some background. We describe various elements of the braid group, including full twists, Young-Jucys-Murphy elements, shuffle braids, and shuffle
twists. In \S\ref{subsec-rouquier} we briefly recall Soergel's categorification of the Hecke algebra and Rouquier's categorification of the braid group. In
\S\ref{subsec-convolutionbackground} we define convolutions of complexes, and give the crucial argument involving the degeneration of a spectral sequence thanks to parity considerations.
In \S\ref{subsec-Hoga-pp} we recall the main result of \cite{Hoga-pp}, a complex which categorifies a renormalized Jones-Wenzl projector, and state its properties. The specifics of the
Soergel-Rouquier construction need not concern the reader, as all we will use in this paper are facts about the braid group and the results of \cite{Hoga-pp}.

In \S\ref{sec-resolve1} we find a convolution description of the full twist $\FT_n$ in terms of certain complexes $D_v$ associated to $v \in \{0,1\}^n$.

In \S\ref{sec-mainResult} we switch to a Hochschild frame of mind. Since Hochschild cohomology of a complex is invariant under conjugation by Rouquier complexes of braids, we will allow
ourselves to freely conjugate complexes. In \S\ref{subsec-markov} we discuss another result of \cite{Hoga-pp} which is an analog of the Markov move on braid closures: a relationship
between the Hochschild homologies of the Jones-Wenzl projector on $n$ strands and the projector on $n-1$ strands. We also discuss reduced complexes in \S\ref{subsec-reducedcomplexes},
finally describing complexes $C_v$ in \S\ref{subsec-ourcomplexes} which may be thought of as reduced versions of conjugates of $D_v$. Finally, in \S\ref{subsec-newRecursion} and
\S\ref{subsec:mainResult} we state and prove the main result, which is a convolution description of the Hochschild cohomology of $C_v$ in terms of the Hochschild cohomologies of
smaller $C_w$, which is an analogue of the recursion of Proposition \ref{prop:fv-intro}.

In \S\ref{sec-numerology} we prove some combinatorial results which justify Theorem \ref{thm-degreezerosolve-intro}, our closed form solution for $\HHH^0(\FT)$, and show that the two
recursive formulas agree.

In the appendix, we include without proof some computations for other $(n,m)$ torus links. These were obtained by techniques entirely analogous to the computation for $(n,n)$ torus links, and many of them have not appeared in the literature before.

{\bf Acknowledgments} The authors would like to thank Eugene Gorsky, Alexei Oblomkov, Pavel Etingof, and Andrei Negut for enlightening conversations. A substantial amount of this work was
completed during the second author's visit to the University of Oregon during the summer of 2015; we are indebted to the UO math department for its hospitality and support, and
apologize for filling the lounge chalkboard with half twists. Both authors would like to thank Yeppie for her excellent sandwiches, which kept us full of hope and inspiration.

\subsection{Notation}
\label{subsec:notation}

We collect here some of our notational conventions, for the reader's convenience.  Unfamiliar concepts will be explained in due course.  Soergel bimodules are graded.  We denote by $(1)$ the grading shift, so that $M(1)^i = M^{i+1}$.  We let $Q=(-1)$ denote the functor which \emph{increases} the degree of each element.  Complexes of Soergel bimodules are bigraded.  The differentials always preserve the bimodule degree, and increase homological degree by 1.  The shift in homological degree is denoted by $\ip{1}$, so that $C(a)\ip{b}^{i,j} = C^{i+a,j+b}$.  We denote by $\K^b(\AC)$ the homotopy category of finite complexes over an additive category $\AC$.  Isomorphism in $\K^b(\AC)$, that is, chain homotopy equivalence, is denoted by $\simeq$.  The existence of a distinguished triangle
\[
A\rightarrow B\rightarrow C\buildrel \d\over \rightarrow A\ip{1}
\]
will be indicated by writing $B\simeq (C \buildrel \d\over \rightarrow A)$.  We also let $T = \ip{-1}$ denote the functor which increases homological degree by 1.  If $\b$ is a braid, we denote the braid exponent by $e(\b)$; this is the signed number of crossings in a diagram representing $\b$.  The Rouquier complex $F(\b)$ is normalized so that if $\b$ is a positive braid, then there is a chain map $(TQ\inv)^{e(\b)}  R\rightarrow F(\b)$ which is the inclusion of the degree $e(\b)$ chain bimodule, whereas if $\b$ is a negative braid, there is a chain map $F(\b)\rightarrow (TQ\inv)^{e(\b)}R$ which is the projection onto the degree $e(\b)$ bimodule.  Note, in \cite{AbHog-pp} and \cite{Hoga-pp}, the shifts $(k)$ and $\ip{\ell}$ would have been denoted $(-k)$ and $\ip{-\ell}$, respectively.

Hochschild cohomology gives rise to a functor $\HH$ whose input is a graded bimodule, and whose output is a bigraded vector space.  The additional grading is called the Hochschild grading, and shifts in the Hochschild grading are denoted by $A$.   Extending to complexes gives a functor from complexes of graded bimodules to complexes of bigraded vector spaces.   These are triply graded objects, so all together we have the shift functors $Q,A,T$.  If $C$ is a complex of bimodules, then the homology of $\HH(C)$ is denoted by $\HHH(C)$.

We also find it convenient to introduce $t=T^2Q^{-2}$, $q=Q^2$, and $a=AQ^{-2}$.  One might call these the \emph{geometric variables}, since they appear most naturally in the connection with Hilbert schemes.  When convenient, we express our degree shifts and Poincar\'e series in terms of these variables.



%% file: FTHHHbackground.tex
\section{Background and key tools} \label{sec-background}

\subsection{Braids} \label{subsec-braids}

Let $\Br_n$ denote the braid group with $n$ strands. The generators will be denoted by $\s_i$, for $1 \le i \le n-1$, and drawn as an \emph{overcrossing} of the $i$-th and $(i+1)$-st
strands. The overcrossing $\s_i$ and its inverse, the \emph{undercrossing} $\s_i\inv$, are depicted below.
\[
\sigma_i \ \ = \ \ \begin{minipage}{1.4in}
\labellist
\small
\pinlabel $i-1$ at -12 12
\pinlabel $n-i-1$ at 62 12
\endlabellist
\begin{center}\includegraphics[scale=1]{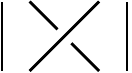}\end{center} 
\end{minipage}
\hskip1in
\sigma_i\inv \ \ = \ \ \begin{minipage}{1.4in}
\labellist
\small
\pinlabel $i-1$ at -12 12
\pinlabel $n-i-1$ at 62 12
\endlabellist
\begin{center}\includegraphics[scale=1]{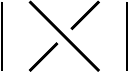}\end{center} 
\end{minipage}
\]
A labelled strand denotes the corresponding number of parallel copies of that strand.  We will always draw our braids in a rectangle, with $n$ boundary points on the top and bottom. Composition of braids is given by vertical stacking, so that $\b\b'$ is $\b$ on top of
$\b'$. There is a group homomorphism $e:\Br_n \to \ZM$ sending $\s_i^\pm \mapsto \pm 1$. The integer $e(\b)$ is called the \emph{braid exponent} of $\b$.

A braid is \emph{positive} if it has an expression only involving overcrossings, and \emph{negative} if it has an expression only involving undercrossings. Given an element $w$ of the
symmetric group $S_n$, its \emph{positive braid lift} in $\Br_n$ is the product $\s_{i_1} \s_{i_2} \cdots \s_{i_d}$, where $s_{i_1} s_{i_2} \cdots s_{i_d}$ is a reduced expression for
$w$ in terms of the usual Coxeter generators $\{s_i\}$ of $S_n$. This element is independent of the choice of reduced expression. Its \emph{negative braid lift} is $\s_{i_1}\inv \cdots \s_{i_d}\inv$.

\begin{defn} \label{def-symmetries-braid}
We define the following symmetries of $\Br_n$:
\begin{enumerate}
\item Rotation about the vertical axis: Let $\t \co \Br_n \to \Br_n$ satisfy $\t(\s_i) = \s_{n-i}$ and $\t(\a \b) = \t(\a) \t(\b)$. Then $\t$ is an involution.
\item Rotation about the horizontal axis: Let $\w \co \Br_n \to \Br_n$ satisfy $\w(\s_i) = \s_i$ and $\w(\a \b) = \w(\b) \w(\a)$. Then $\w$ is an antiinvolution.
\item Reflection across a horizontal plane: Let $(-)^\vee \co \Br_n \to \Br_n$ satisfy $\s_i^{\vee} = \s_i\inv$ and $(\a \b)^\vee = \b^\vee \a^\vee$. Then $(-)^\vee$ is an antiinvolution, and is just another notation for taking the inverse braid.
\item Crossing swap: Let $(-)^L \co \Br_n \to \Br_n$ satisfy $\s_i^L = \s_i\inv$ and $(\a \b)^L = (\a)^L (\b)^L$. Then $(-)^L$ is an involution, and $(\b)^L = \w(\b)^\vee$.
\end{enumerate}
\end{defn}

The letter $L$ indicates that $\b^L$ this is the \emph{left-handed version} of the braid $\b$. This swaps the positive and negative braid lifts of an element of $S_n$. Note that $\t$ and
$\w$ preserve positive braids, while $(-)^\vee$ and $(-)^L$ swap positive braids and negative braids. These symmetries all commute with each other.

We let $\sqcup \co \Br_k \times \Br_l \to \Br_{k+l}$ denote the homomorphism given by horizontal concatenation.

\subsection{Shuffle braids} \label{subsec-shuffles}

\begin{defn} \label{def-shuffle} A \emph{shuffle permutation} is a permutation $\pi \in S_n$ which is a minimal length coset representative for some coset in $S_n / (S_k \times S_\ell)$,
for some $0 \le k,\ell \le n$ with $k+\ell = n$.  Said differently, a shuffle permutation preserves the ordering of $\{1,\ldots,k\}$ and $\{k+1,\cdots,n\}$ for some $k$, but ``shuffles'' these two sets together.

Let $v \in \{0,1\}^n$ be a sequence with $k$ zeroes and $\ell$ ones. We call $v$ a \emph{shuffle}. There is a corresponding shuffle permutation $\pi_v$, a minimal coset representation
for $S_n / (S_k \times S_\ell)$, for which $\pi_v(\{1, \ldots, k\})$ gives the locations of the zeroes, and $\pi_v(\{k+1, \ldots, n\})$ gives the locations of the ones.
\end{defn}

Note that a shuffle permutation can come from a shuffle in multiple different ways. For example, the identity element is a minimal coset representative for $S_k \times S_\ell$ for every
$k$ and $l$ with $k+\ell=n$; whenever all the zeroes come before all the ones, $\pi_{0 \cdots 01 \cdots 1}$ is the identity. When the shuffle $v \in \{0,1\}^n$ is understood, $k$
will always refer to the number of zeroes, and $\ell$ to the number of ones.

\begin{ex} The shuffle permutation $\pi_{1\cdots 1 0}$ is the $n$-cycle $(n,n-1,\ldots,2,1)$. \end{ex}

\begin{defn}\label{def-shuffleBraids}
For each $v\in\{0,1\}^n$, let $\b_v$ denote the positive braid lift of $\pi_v$. Let $\Tw_v:=\omega(\b_v)\b_v$, the \emph{shuffle twist}, denote the positive pure braid obtained by gluing $\b_v$ with its rotation. \end{defn}

\begin{ex}
If $v=(0101100)$, then the shuffle permutation $\pi_v$, its positive braid lift, and the associated pure braid are pictured as:
\[
\pi_v = \ig{1}{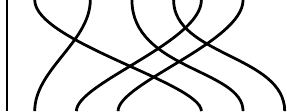} \hskip.5in \b_v = \ig{1}{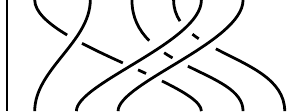} \hskip.5in \Tw_v = \ig{1}{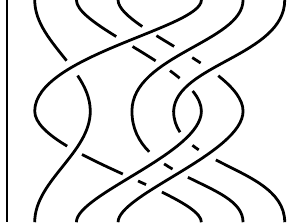}.
\]
Note that $v$ partitions the strands in these diagrams into two subsets: the 0-strands and the 1-strands. In $\b_v$ the 0-strands cross over the 1-strands. In the $\w(\b_v)$ portion of $\Tw_v$, they cross back under.
\end{ex}

The following gives a useful recursive description of the braids $\Tw_v$:

\begin{prop}\label{prop-shuffleRecursion}
Let $v$ be a shuffle with $k$ zeroes and $\ell$ ones, with $k+\ell = n$. Let $\cdot$ denote concatenation of shuffles, so that $v \cdot 0$ and $v \cdot 1$ are the two shuffles of length $n+1$ which extend $v$. Then 
\[
\Tw_{v\cdot 0} = 
\begin{minipage}{1.1in}
\labellist
\small
\pinlabel $\Tw_v$ at 18 25
\pinlabel $k$ at -5 8
\pinlabel $1$ at 10 8
\pinlabel $\ell$ at 35 0
\endlabellist
\begin{center}\includegraphics[scale=1]{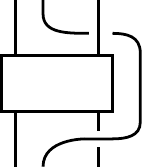}\end{center} 
\end{minipage}
\hskip1in
\Tw_{w\cdot 1} = 
\begin{minipage}{1.1in}
\labellist
\small
\pinlabel $\Tw_v$ at 18 25
\pinlabel $k$ at 0 8
\pinlabel $\ell$ at 18 8
\pinlabel $1$ at 32 8
\endlabellist
\begin{center}\includegraphics[scale=1]{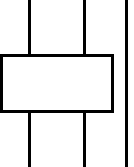}.\end{center}
\end{minipage}
\]
Here a strand labeled $\ell$ actually represents $\ell$ strands cabled together in the usual way, so that the ``thick crossings" in the expression for $\Tw_{v \cdot 0}$ each represent $\ell$ ordinary crossings. This recursion, together with the base cases $\Tw_{0}=\Tw_{1}=1$, produces all the shuffle twists $\Tw_v$.  There is a similar such recursion which describes $\Tw_{0\cdot v}$ and $\Tw_{1\cdot v}$.
\end{prop}

\begin{proof}
Graphically obvious.
\end{proof}

We now discuss the behavior of the shuffle braids with respect to the symmetries of the braid group.

\begin{defn}\label{def-shuffleSymmetries}
For each $v\in\{0,1\}^n$, let $r(v)$ denote the sequence obtained by reversing the order, so that $r(v)_i = v_{n+1-i}$.  Let $v^\ast$ be the sequence obtained by swapping the 1's with 0's and vice versa: $(v^\ast)_i = 1-v_i$.
\end{defn}

\begin{prop}\label{prop-shuffleSymmetries}
Let $v\in\{0,1\}^n$ be given.  Then
\begin{enumerate}
\item $\tau(\b_v) = \b_{r(v)^\ast}$
\item $\omega(\b_v\inv) = \b_v^L$
\end{enumerate}
\end{prop}
\begin{proof}
Clear.
\end{proof}

\subsection{Half twists and shuffle braids} \label{subsec-HTshuffles}

Let $\HT=\HT_n\in \Br_n$ denote the half twist braid
\[
\HT_n = \sigma_1(\sigma_2\sigma_1)(\sigma_3\sigma_2\sigma_1)\cdots(\sigma_{n-1}\cdots \sigma_2\sigma_1).
\] 
The full twist is $\FT_n = \HT_n^2$, and is central in the braid group.  This implies that the mapping $\b\mapsto \HT \b \HT\inv$ defines an involution on the braid group.  Indeed, $\HT \b\HT\inv = \tau(\b)$.  We leave the proof of this fact to this reader, as it is elementary.  It is also elementary that $\HT_n$ is fixed by $\tau$ and $\omega$.

We will need to know how the shuffle braids interact with $\HT$ and $\FT$:

\begin{prop}\label{prop-HTshuffle}
Let $v\in \{0,1\}^n$ be given.  Let $k$ and $\ell$ be the number of zeroes and ones in $v$, respectively.  Then
\begin{enumerate}
\item $\HT_n \gamma_v= \gamma_{r(v)^\ast}\HT_n$, where $\gamma_v$ is any of the braids $\b_v$, $\w(\b_v)$, $\b_v^L$, or $\w(\b_v)^L$.
\item $\HT_n \b_v^L = \b_{r(v)}(\HT_k\sqcup \HT_{\ell})$.
\item $\FT_n \Tw_v^L \sim \Tw_{r(v)}(\FT_k\sqcup \FT_{\ell})$.
\end{enumerate}
Here $\sim$ denotes that the given braids are equivalent modulo conjugation.
\end{prop}
\begin{proof}
Statement (1) follows from Proposition \ref{prop-shuffleSymmetries}, since conjugation by $\HT_n$ acts on the braid group by $\t$ (180 degree rotation about a vertical axis).  

The idea of the proof of (2) is best illustrated with an example.  For instance, when $v=(1100101)$ we have
\[
\begin{minipage}{1.2in}
\labellist
\small
\pinlabel $\HT_{k+\ell}$ at 45 40
\endlabellist
\begin{center}\includegraphics[scale=.8]{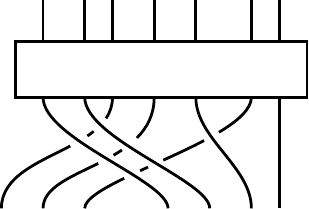}\end{center} 
\end{minipage}
\simeq 
\begin{minipage}{1.2in}
\labellist
\small
\pinlabel $\HT_{k+\ell}$ at 45 20
\endlabellist
\begin{center}\includegraphics[scale=.8]{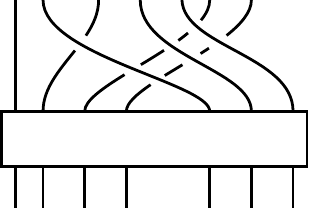}\end{center} 
\end{minipage}
\simeq
\begin{minipage}{1.4in}
\labellist
\small
\pinlabel $\HT_{k}$ at 17 16
\pinlabel $\HT_{\ell}$ at 70 16
\endlabellist
\begin{center}\includegraphics[scale=.8]{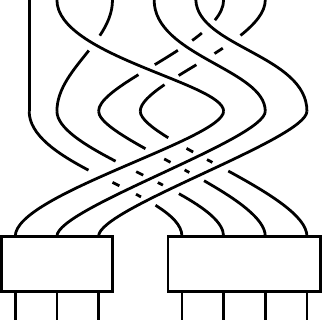}\end{center} 
\end{minipage}
\simeq
\begin{minipage}{1.4in}
\labellist
\small
\pinlabel $\HT_{k}$ at 17 16
\pinlabel $\HT_{\ell}$ at 70 16
\endlabellist
\begin{center}\includegraphics[scale=.8]{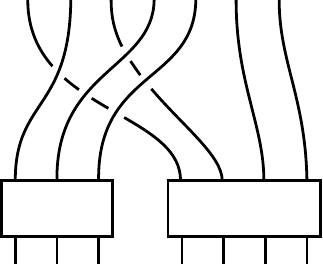}\end{center} 
\end{minipage}.
\]
In the second diagram, the left-handed shuffle braid corresponds to $r(v)^\ast = (0101100)$ by statement (1) of the proposition.  In the third diagram we have simply rewritten the half twist in terms of the ``thick crossing'' between $4$ cabled strands and 3 cabled strands.  This is a well-known identity in the braid group.  In the last diagram we have performed an isotopy.

Finally, statement (3) follows from (2).  First, note that an application of $\omega$ to statement (2) yields
\[
\omega(\b_{v}^L)\HT = (\HT_k\sqcup \HT_\ell)\omega(\b_{r(v)})
\]
Then observe:
\begin{eqnarray*}
\FT \Tw^L_v &=& \FT \omega(\b_v^L)\b_v^L\\
&=&  \omega(\b_{v}^L)\FT \b_{v}^L \\
&=&  \omega(\b_{v}^L)\HT \HT \b_{v}^L \\
&=&  (\HT_k\sqcup \HT_\ell)\omega(\b_{r(v)}) \b_{r(v)} (\HT_k\sqcup \HT_\ell)\\
&\sim &  \omega(\b_{r(v)}) \b_{r(v)} (\FT_k\sqcup \FT_\ell) \\
\end{eqnarray*}
The first equality holds since $\FT_n\in \Br_n$ is central.  The second holds since $\FT = \HT \HT$.  The third holds by (1) and (2).  Finally the last $\sim$ holds by transferring the $(\HT_k\sqcup \HT_\ell)$ to the right-hand side (recall that $\b\sim \b'$ means $\b$ is conjugate to $\b'$).  This completes the proof.
\end{proof}

\subsection{Rouquier complexes} \label{subsec-rouquier}

Let $R = R_n = \QM[x_1, \ldots, x_n]$ be the polynomial ring in $n$ variables, graded so that $\deg x_i = 2$. This is the polynomial ring associated to the standard $n$-dimensional
representation of $S_n$ over $\QM$. Given a graded $R$-bimodule $M$, we let $M(1)$ denote the shifted bimodule for which $M(1)^d = M^{1+d}$, where $M^d$ denotes the degree $d$ part of
$M$. We denote tensor product of graded bimodules over $R$ simply by juxtaposition: $M \otimes_R N = MN$. We often let $\1 = \1_n$ denote the bimodule $R$, which is the monoidal identity.

For each $i$ with $1 \le i \le n-1$, let $B_i$ denote the graded $R$-bimodule \[B_i := R \ot_{R^i} R (1), \] where $R^i$ denotes the subring of polynomials invariant under the reflection
$s_i = (i,i+1)$. Let $\SBim_n$ denote the category of Soergel bimodules associated to $R$. This is the full graded monoidal additive Karoubian subcategory of graded $R$-bimodules
generated by $B_i$ for $1 \le i \le n-1$. Thus, its objects are direct sums of grading shifts of direct summands of tensor products of $B_i$. Let $\CC^b(\SBim)$ denote the category of bounded complexes of Soergel bimodules modulo homotopy.

Associated to each braid word $\b$ we have the Rouquier complex $F(\b)$ in $\CC^b(\SBim)$, defined by
\[
F(\s_i) = (\un{B_i}\to R(1)) \ \ \ \ \ \ \ \ \ \ \ \ \ F(\s_i\inv)=(R(-1)\to \un{B_i})  
\]
together with $F(\b\b')=F(\b) F(\b')$. The underline indicates which object lies in homological degree $0$. Rouquier proved that there is a canonical homotopy equivalence between $F(\b)$ and $F(\b')$ when $\b$ and $\b'$ are braid words expressing the same braid. A more direct proof which works over $\ZM$ can also be found in \cite{EKra}.

\begin{remark} Recall the notation $Q,T$ for gradings shifts (see \S \ref{subsec:notation}). One reason why $t$ is more natural that $T$ is that any Rouquier complex $F(\b)$ always has a unique copy of $R$ which appears in homological degree $e(\b)$ and internal degree $-e(\b)$, so that this copy of $R$ appears with shift $t^{\frac{1}{2}e(\b)}$, where $e(\b)$ is the braid exponent.  It was proven in \cite{EWHodge} that Rouquier complexes for reduced expressions are \emph{perverse} (when one works in characteristic zero).  A complex is perverse if each indecomposable bimodule in the complex appears with a grading shift equal to its homological degree, or equivalently, that the grading and homological shifts are described only as powers of $(TQ\inv)$. Note that the Rouquier complex for the full twist is not perverse, nor are shuffle twists. We will not use any perversity results in this paper.  Nonetheless, we will express our shifts using the variables $t=T^2Q^{-2}$ and $q=Q^2$.\end{remark}

The symmetries of the braid group lift to symmetries of $\SBim_n$ and its homotopy category. Let $\t \co R \to R$ denote the map sending $x_i \mapsto x_{n+1-i}$.

\begin{defn}\label{def-symmetries-rouquier}
We define the following symmetries of $\SBim_n$:
\begin{enumerate}
\item Rotation about vertical axis: let $\t:\SBim_n \to \SBim_n$ denote the covariant graded monoidal functor induced by the Dynkin automorphism of $S_n$.  That is, $\t: B_i \mapsto B_{n-i}$, and satisfies $\t(MN) = \t(M)\t(N)$ and $\t(M(1)) = \t(M)(1)$.
\item Rotation about horizontal axis: let $\w: \SBim_n \to \SBim_n$ denote the covariant graded anti-monoidal functor which sends $B_i \mapsto B_i$ and satisfies $\w(MN)=\w(N)\w(M)$ and $\w(M(1)) = \w(M)(1)$.
\item Reflection across a horizontal plane: let $(-)^\vee : \SBim_n \to \SBim_n^{\op}$ denote the contravariant anti-graded anti-monoidal ``duality" functor on $\SBim_n$, which sends $B_i \mapsto B_i$ and satisfies $(MN)^\vee = N^\vee M^\vee$ and $M(1)^\vee = M^\vee(-1)$.
\item Crossing swap: Let $(-)^L \co \SBim_n \to \SBim_n$ denote the contravariant anti-graded monoidal functor $(-)^\vee \circ \w$.
\end{enumerate}
\end{defn}

These symmetries commute up to canonical isomorphism.

\begin{prop} These symmetries intertwine the braid symmetries from Definition \ref{def-symmetries-braid}, under Rouquier's map $F$, via a canonical isomorphism. In particular, $F(\b)^\vee \cong F(\b\inv)$.
\end{prop}

\begin{proof} This is easy. Although we have not stated explicitly what the differentials in $F(\s_i)$ and $F(\s_i\inv)$ are, they live in one-dimensional morphism spaces, and are
interchanged by duality. \end{proof}

There is an isomorphism of rings $R_k \boxtimes R_l \to R_{k+l}$ given by renaming the variables, where $\boxtimes$ denotes tensor product over $\QM$. Correspondingly, there is an
inclusion functor $\sqcup \co \SBim_k \boxtimes \SBim_l \to \SBim_{k+l}$, which sends $B_i \boxtimes \1_l$ to $B_i$ and sends $\1_k \boxtimes B_i$ to $B_{k+i}$. This intertwines with the
map $\sqcup \co \Br_k \times \Br_l \to \Br_{k+l}$ after applying Rouquier's map $F$.

It was proven by Soergel that morphisms between objects in $\SBim_n$ are free as left or right modules over $R_n$. Using Soergel's Hom formula \cite[Thm 5.15]{Soer07}, one can prove that $\sqcup$ is actually fully faithful, after identifying $R_k \boxtimes R_l$ with $R_n$. Another way of phrasing this result is that the inclusion
$\SBim_i \to \SBim_n$ for $i < n$, which comes from the functor $(-) \sqcup \1_{n-i}$, is fully faithful after base change along the inclusion from $R_i$ to $R_n$. See \cite[Remark 3.19]{EWGr4sb} for further discussion.

\subsection{Complexes and convolutions}\label{subsec-convolutionbackground}

We may write $\ip{1}$ for the homological shift of a complex, so that the homological degree $d$ part of $F\ip{1}$ is the homological degree $d+1$ part of $F$.  By convention, $\ip{1}$ also negates the differential.

We now introduce some notation which we will be used exhaustively throughout.  To motivate it, we begin with an example.  Suppose $A$ and $B$ are complexes, and $f:A\rightarrow B$ is a chain map.  The mapping cone $C_f$ is the chain complex $(C_f)_k = A_{k+1}\oplus B_k$ with differential given by the matrix $\left[\begin{smallmatrix}-d_A &0\\f&d_B\end{smallmatrix}\right]$.  In other words $C_f = A\ip{1}\oplus B$ with an additional component of the differential from $A\ip{1}$ to $B$, given by $f$.  Note that the additional sign on the differential of $A\ip{1}$ is necessary for the differential on $C_f$ to satisfy $d^2=0$.   We prefer to keep track of the homological degree shift explicitly, so that the mapping cone can be written as
\[
C_f = (A\ip{1}\buildrel f\over\longrightarrow B).
\]
We will also say that $C_f = A\ip{1}\oplus B$ with twisted differential.  This notation will come in handy when we later consider mapping cones of mapping cones, and so on.  

For instance, this notation allows us to use explicit shifts instead of
underlines in a complex, so that we may write
\[ F(\s_i) = (B_i \to R(1)\ip{-1}) \ \ \ \ \ \ \ \ \ \ \ \ \ F(\s_i\inv)=(R(-1)\ip{1} \to B_i). \]

The general way to describe an iterated cone is using the idea of a convolution of complexes. Let $F_j$ ($j \in J$) be complexes of $R$-bimodules indexed by a finite partially-ordered
set $J$. Let $d_j$ denote the differential on the complex $F_j$ (which, in our notation, is a map of bigraded $R$-bimodules of homological degree $+1$ and graded degree $0$). Let $E =
\oplus_{j \in J} F_j$ be a bigraded $R$-bimodule, and let $d$ be a differential on $E$ such that
\begin{itemize} \item restricted to a map $F_j \to F_j$, $d$ agrees with $d_j$, and \item restricted to a map $F_j \to F_{j'}$, $d$ is zero unless $j \le j'$. \end{itemize} Then $E$ is
called a \emph{convolution} of the complexes $F_j$, as is any complex which is homotopy equivalent to $E$. We may write $d = \sum_{i\leq j} d_{ji}$, where $d_{ji}$ is the component of the differential mapping $F_i$ to $F_j$.

We refer to $F_j$ as the \emph{subquotients} of the convolution $E$. We say that $E = \oplus_{j \in J} F_j$ with \emph{twisted differential}, indicating that the differential is not just the direct sum of the differentials on each summands. We say that this twisted differential respects the partial order on $J$ because $d_{ji}=0$ for $i \nleq j$.

\begin{remark}\label{rmk:cones}
By abuse of language, we will refer to a two term convolution $E = (A \buildrel f\over \longrightarrow  B)$ as a mapping cone.  Note that, strictly speaking $f$ is not a chain map from $A$ to $B$, but rather a chain map $A\ip{-1}\rightarrow B$.  Here $J$ has two elements, with the order determined by the arrow.  A general convolution can be described as an iterated cone of complexes, where each $F_j$ is added one at a time.
\end{remark}

In practice, one can often show indirectly that a complex $E$ is a convolution of other complexes $F_j$, in which case the components $d_{ji}$ of the differential may be difficult to write down for $i \ne j$ (the task is complicated further by the presence of homotopy equivalences). In particular, this makes it difficult to compute the homology $H(E)$. Thankfully, a parity argument will come to the rescue in this paper.

\begin{proposition}\label{prop:evenConvolutions}
Suppose $E = \bigoplus_{j\in J} F_j$ with twisted differential, for some finite partially ordered set $J$.  Suppose the homology $H(F_j)$ is supported in even homological degrees, for all $j\in J$.  Then $H(E)\cong \bigoplus_j H(F_j)$.
\end{proposition}

\begin{proof}
We induct on the cardinality of $J$.  In the base case $J=\{j\}$, we have $E=F_j$, and the statement is trivial.  Now, assume by induction that we have proved the result for partially ordered sets of cardinality $r$, and let $J$ be a partially ordered set of cardinality $r+1$. Let $j\in J$ be maximal.  Set $B:=F_j$ and $A=\bigoplus_{i\in J\smallsetminus \{j\}} F_j$ with twisted differential.  Note that $E = A\oplus B$ with twisted differential:
\[
E = (A\buildrel \d\over \longrightarrow B)
\]
for some map $\d$ of homological degree $+1$. We have $H(A) \cong \bigoplus_{i\neq j} H(F_i)$ by induction, so we must prove that $H(E)\cong H(A)\oplus H(B)$.

The short exact sequence $0\rightarrow B\rightarrow E \rightarrow A\rightarrow 0$ gives rise to a long exact sequence
\[
\cdots \rightarrow H^{k-1}(A)\rightarrow H^{k}(B)\rightarrow H^{k}(E)\rightarrow H^k(A)\rightarrow H^{k+1}(B)\rightarrow \cdots.
\]
Our parity assumption implies that $H^{k}(A)=H^{k}(B)=0$ when $k$ is odd. Thus $H^k(E)=0$ when $k$ is odd. When $k$ is even we have a short exact sequence
\[
0\rightarrow H^{k}(B)\rightarrow H^{k}(E)\rightarrow H^{k}(A)\rightarrow 0.
\]
If we work over a field, then this short exact sequence splits.  This completes the inductive step, and completes the proof.
\end{proof}

\begin{remark}
In general, there is a spectral sequence converging to $H(E)$, whose $E_2$ page is $\bigoplus_j H(F_j)$.  If $H(F_i)$ is even, then the subsequent differentials (which have odd homological degree) must all vanish.  This gives an alternate proof of the above.
\end{remark}

\begin{remark}
The above presents a ``computation-free and serendipitous'' approach to computing homology groups. Suppose we wish to compute the homology of a chain complex $E$.  We may get lucky and discover a filtration on $E$ whose successive quotients are supported in even homological degrees.  In this case, Proposition \ref{prop:evenConvolutions} says that $H(E)$ simply splits as a direct sum of these homology groups.  In this paper we are extraordinarily lucky in this regard.
\end{remark}

\subsection{Categorified symmetrizers}\label{subsec-Hoga-pp}
In this subsection we recall the constructions of the second author in \cite{Hoga-pp}, and extract from them a finite complex $K_n\in \K^b(\SBim_n)$ which will play an essential role in this paper.  First, define the following complexes:
\begin{defn}\label{def-XY}
Let $X=X_n=F(\sigma_{n-1}\cdots \sigma_1)$ and $Y=Y_n=F(\sigma_{n-1}\inv \cdots \sigma_1\inv)$ denote the Rouquier complexes associated to the positive and negative braid lifts of the standard $n$-cycle $(n,n-1,\ldots,2,1)$. In other words, $X = F(\b_v)$ and $Y = F(\b_v^L)$ for $v = (11 \cdots 10)$. Note that $\t(X)=\w(X)=Y\inv$, $\t(Y)=\w(Y)=X\inv$.
\end{defn}

In this paper we adopt a graphical notation for certain complexes of Soergel bimodules.  We will denote a braid and its Rouquier complex similarly, so for example pictures such as 
\[
X_4 = \ig{1}{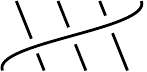} \ \ \ \ \ \ \text{ and } \ \ \ \ \ \  Y_4 = \ig{1}{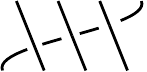}
\]
will denote the complexes $X_n$ and $Y_n$ of Definition \ref{def-XY}.  The tensor product of complexes corresponds to vertical stacking.  For example we have
\begin{equation}\label{eqn:yjm4}
X_4Y_4\inv = \ig{1}{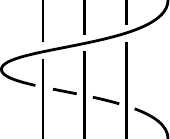} \ \in \Br_4
\end{equation}
Complexes of the form $XY\inv$ play a very special role in this paper.  They are the Rouquier complexes associated to the Young-Jucys-Murphy braids.  Note that braids corresponding to $X_kY_k\inv$ generate a commutative subgroup of the braid group, and their product is the full twist.  The following defines a family of complexes $K_n$ which are compatible with these braids, in a particular sense.

\begin{proposition}\label{prop:K}
There exists a family of finite complexes $K_n \in \CC^b(\SBim_n)$ ($n\geq 1$) such that:
\begin{enumerate}
\item $K_1 = R$.
\item We have $K_{n-1}XY\inv \simeq (K_n\rightarrow q K_{n-1})$.  Graphically this is
\begin{equation} \label{eqn-YJM-cone}
\begin{minipage}{1in}
\labellist
\small
\pinlabel $K_{n-1}$ at 27 39
\endlabellist
\includegraphics[scale=1]{diagrams/KwithYJM0}
\end{minipage}
\ \ \simeq \ \ \left( \
\begin{minipage}{1in}
\labellist
\small
\pinlabel $K_{n}$ at 28 29
\endlabellist
\includegraphics[scale=1]{diagrams/KwithYJM1}
\end{minipage}
\longrightarrow
q \ \begin{minipage}{1in}
\labellist
\small
\pinlabel $K_{n-1}$ at 28 29
\endlabellist
\includegraphics[scale=1]{diagrams/KwithYJM2}
\end{minipage}
\right)
\end{equation}
\item $K_n$ kills all generating Bott-Samelson bimodules in $\SBim_n$.  That is, $K_n B_i \simeq 0 \simeq B_i K_n$ for $1 \le i \le n-1$.
\end{enumerate}
\end{proposition}

Note that finiteness of the $K_n$ follows from the recursion (2).  

\begin{proof}
In \cite{Hoga-pp} the second author defined defined a family of complexes $P_n\in \K^-(\SBim_n)$ such that
\begin{itemize}
\item $P_n$ kills all Bott-Samelsons: $P_n\otimes B_i\simeq 0\simeq B_i \otimes P_n$ for all $1\leq i\leq n-1$.
\item Any other complex $M\in \K^-(\SBim_n)$ kills Bott-Samelsons if and only if $P_n\otimes M\simeq M\simeq P_n\otimes M$.
\end{itemize}
Further, $P_n$ is uniquely characterized up to homotopy equivalence by these properties.  The complexes $P_n$ can be thought of as categorical projections onto the sign representation of the Hecke algebra:

We will construct the complexes $K_n$ inductively.  First, set $K_1 = R_1$.  Now, assume $K_{n-1}$ has been constructed for $n\geq 2$.   In \S 4 of \cite{Hoga-pp} it was shown that there is a chain map $\psi: q P_{n-1}\rightarrow P_{n-1}XY\inv$ such that the mapping cone $Q_n:=\Cone(\psi)$ kills Bott-Samelsons\footnote{Actually, in \cite{Hoga-pp} the renormalized Rouquier complexes $F'(\b):=t^{-\frac{1}{2}e(\b)}F(\b)$ are used, where $e(\b)$ is the braid exponent.  To match the conventions, set $X'=t^{\frac{1}{2}(1-n)}X$ and $Y'=t^{\frac{1}{2}(n-1)}Y$.  Then the chain map constructed in \cite{Hoga-pp} is actually $\psi':qt^{1-n}P_{n-1}\rightarrow P_{n-1}X'(Y')\inv$.  Clearly this gives rise to $\psi:q P_{n-1}\rightarrow P_{n-1}XY\inv$ as claimed.}.  Graphically, this is:
\begin{equation} \label{eqn-Qn-cone}
\begin{minipage}{1in}
\labellist
\small
\pinlabel $Q_{n}$ at 28 29
\endlabellist
\includegraphics[scale=1]{diagrams/KwithYJM1}
\end{minipage}
\ \ \simeq \ \ 
 \left(
q \ \begin{minipage}{.85in}
\labellist
\small
\pinlabel $P_{n-1}$ at 28 29
\endlabellist
\includegraphics[scale=1]{diagrams/KwithYJM2}
\end{minipage}\ip{1}
\ \buildrel \psi\over\longrightarrow \ 
\begin{minipage}{1in}
\labellist
\small
\pinlabel $P_{n-1}$ at 27 39
\endlabellist
\includegraphics[scale=1]{diagrams/KwithYJM0}
\end{minipage}
\right).
\end{equation}

From the above characterization of $P_{n-1}$, the fact that $K_{n-1}$ kills the Bott-Samelsons in $\SBim_{n-1}$ implies that $K_{n-1}P_{n-1}\simeq K_{n-1}$.  We define $K_n$ by tensoring Equation (\ref{eqn-Qn-cone}) on the left with $K_{n-1}$ and applying the equivalence $K_{n-1}P_{n-1}\simeq K_{n-1}$.  Note that $K_n \simeq K_{n-1}Q_n$.  The recursion (2) is satisfied after rotating triangles.   Clearly $K_n$ kills all Bott-Samelsons on the right since $Q_n$ does.  It was shown in \cite{Hoga-pp} that a complex kills Bott-Samelsons from the right if and only if it kills all Bott-Samelsons from the left.  This gives statement (3).
\end{proof}

\begin{lemma}\label{lemma-RouquierAbsorbing}
The complex $K_n$ absorbs Rouquier complexes: if $\b$ is a braid, then $K_n F(\b) \simeq t^{\frac{1}{2}e(\b)} K_n \simeq F(\b)K_n$. Recall that $e(\b)$ is the braid exponent, which records the number of overcrossings minus the number of undercrossings in a braid word.  
\end{lemma}
\begin{proof}
It suffices to prove the result for $\b=\sigma_i^\pm$.  In this case the claim is obvious since $K_n$ kills $B_i$, hence the only term of $K_n F(\sigma_i^\pm)$ which survives is $K_n(\pm 1)\ip{\mp 1}$.
\end{proof}

\begin{remark}
It is sometimes also useful to consider the following equivalence:
\begin{equation} \label{eqn-Kn-cone}
\begin{minipage}{1in}
\labellist
\small
\pinlabel $K_{n-1}$ at 28 32
\endlabellist
\includegraphics[scale=1]{diagrams/KrecursionRight}
\end{minipage}
\ \ \simeq \ \ 
\left(
t^{\frac{1}{2}(1-n)}\begin{minipage}{1.08in}
\labellist
\small
\pinlabel $K_{n}$ at 27 24
\endlabellist
\includegraphics[scale=1]{diagrams/KrecursionBig}
\end{minipage}
\ \ \longrightarrow \ \ 
q\; \begin{minipage}{1.08in}
\labellist
\small
\pinlabel $K_{n-1}$ at 28 32
\endlabellist
\includegraphics[scale=1]{diagrams/KrecursionLeft}
\end{minipage}
\right).
\end{equation}
This follows from (\ref{eqn-YJM-cone}) by tensoring on the right with $Y$, applying the equivalence $K_n Y\simeq t^{\frac{1}{2}(1-n)}K_n$.
\end{remark}

\begin{ex}\label{ex-K2}
There is a chain map $q F(\sigma_1\inv)\rightarrow F(\sigma_1)$ whose mapping cone is the 4-term complex
\[
K_2 \ \ \  = \ \ \  \begin{diagram} R(-2) &\rTo& \underline{B_1}(-1)&\rTo^{x_2\otimes 1 - 1\otimes x_2}& {B_1}(1)&\rTo & {R}(2).\end{diagram}
\]
There is a projection map $K_2\rightarrow R(-2)\ip{1}$, the mapping cone on which is
\[
(K_2\rightarrow q R) \ \ \simeq \ \ \begin{diagram} \underline{B_1}(-1)&\rTo^{x_2\otimes 1 - 1\otimes x_2}& {B_1}(1)&\rTo& {R}(2).\end{diagram}
\]
By \eqref{eqn-Kn-cone}, this is homotopy equivalent to the full twist $\FT_2 = F(\sigma_1^2)$ on two strands.  This fact is also straightforward to check directly.
\end{ex}

The construction of $K_n$ appears to be asymmetric.  However, $K_n$ is preserved by the symmetries of the Soergel category up to homotopy:
\begin{prop}\label{prop-Ksymmetries}
We have
\[
\w (K_n)\simeq \t(K_n) \simeq K_n.
\]
Further,
\[
t^{n-1} K_n^\vee \simeq  (t^{-\frac{1}{2}}q^{\frac{1}{2}})^{n-1} K_n .
\]
\end{prop}

%% file: FTHHHconvolve1.tex
\section{Resolving the full twist} \label{sec-resolve1}

In this section we give a new expression for the Rouquier complex associated to the full twist braid $\FT_n$.  Our main result is that $\FT_n \in \CC^b(\SBim_n)$ is homotopy equivalent to a convolution whose subquotients are described in terms of shuffle braids and the complexes $K_n$.

\subsection{Young-Jucys-Murphy braids}

First, we describe the Rouquier complexes for Young-Jucys-Murphy braids as convolutions.  We write $y_i=\sigma_{i-1}\cdots \sigma_2\sigma_1\sigma_1\sigma_2\cdots \sigma_{i-1}$ for the $i$-th Young-Jucys-Murphy braid, which is an element of $\Br_i$.  For example, $y_4$ is picture in \eqref{eqn:yjm4}.   We may
also view $y_i$ as an element of $\Br_n$ for any $n \ge i$, which acts on the first $i$ strands; this comes from the inclusion $\Br_i = \Br_i \times 1_{n-i} \subset \Br_i \times \Br_{n-i}
\to \Br_n$.

\begin{defn} \label{defn-E} Let $v \in \{0,1\}^{n}$ be a shuffle on $n$ letters, with $k$ zeroes and $\ell$ ones, so that $k + \ell = n$.  Let $E_{v} \in \CC^b(\SBim)$ be defined as follows:
\begin{equation} \label{eqn-E-defn} E_{v} := F(\b_{v}^L)(\one_k \sqcup K_{\ell})F(\w(\b_{v}^L)). \end{equation} Note that the shuffle braids involved here are left-handed. \end{defn}

\begin{ex} For example, when $v = (1011001)$, the complex $E_{v}$ looks like
 \[
 \begin{minipage}{1.15in}
\labellist
\small
\pinlabel $K_4$ at 65 39
\endlabellist
\begin{center}\includegraphics[scale=.85]{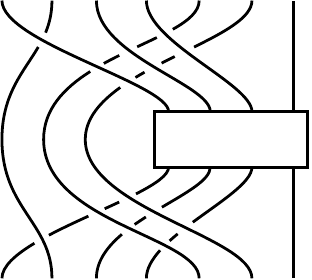}\end{center} 
\end{minipage}.
\] \end{ex}
In this section the only complexes that concern us are those of the form $E_{v\cdot 1}$.  Recall that $\cdot$ denotes concatenation of sequences, so that $v\cdot 1$ ranges over all sequences which end in 1.  Note that (analogously to Proposition \ref{prop-shuffleRecursion}) $\b_{v \cdot 1}$ is equal to $\b_v \sqcup \1_1$, and the right-most strand in $E_{v\cdot 1}$ is a straight vertical line which does not cross over or under any other strands.
\begin{proposition}\label{prop-yjmSimp}
The Rouquier complex $F(y_n)$ satisfies
\[
F(y_n)\simeq \bigoplus_{v\in\{0,1\}^{n-1}}t^{-\binom{\ell}{2}} q^{k}E_{v \cdot 1}
\]
with twisted differential.  As usual, $k$ is the number of zeroes in $v$ and $\ell$ is the number of ones, with $k+\ell=n-1$.  The partial order in this convolution is the antilexicographic order on sequences.
\end{proposition}

\begin{remark}
Consider the ``thick crossing''  between $K_m$ and $n-m$ parallel strands, which one might picture as
\[
\begin{minipage}{.8in}
\labellist
\small
\pinlabel $K_{m}$ at 10 24
\pinlabel $n-m$ at -4 7
\pinlabel $m$ at 36 7
\endlabellist
\begin{center}\includegraphics[scale=1]{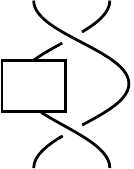}\end{center} 
\end{minipage}.
\]
There is an expression of this complex as a direct sum (with twisted differential) of complexes $E_{v\cdot 1^m}$ with shifts, where $v\in\{0,1\}^{n-m}$.  Proposition \ref{prop-yjmSimp} corresponds to the case $m=1$.  This is the only case that concerns us, so we leave the statement (and proof) for $m>1$ to the reader.
\end{remark}

\begin{proof}
We prove this by induction on $n \ge 1$.   In the base case we have $y_1=1$ and $F(y_1)=K_1$.  There is exactly one element of $\{0,1\}^0$, the empty sequence, so that the sum on the RHS has one term $E_{\emptyset \cdot 1}$, and that term is $K_1$.  This establishes the base case.

Assume by induction that the result holds for $n\geq 1$.  Note that $y_{n+1} = \sigma_n y_n \sigma_n$.  By induction, we have
\begin{equation}\label{eq-yjmInduction}
F(\sigma_n)F(y_n) F(\sigma_n) \simeq \bigoplus_{v\in\{0,1\}^{n-1}} G_v F(\s_n) (E_{v \cdot 1} \sqcup \1_1) F(\sigma_n)
\end{equation}
with twisted differential, for some grading shifts $G_v$.  Each of the above summands can be rewritten as follows: for fixed $v\in\{0,1\}^{n-1}$, let $k$ denote the number of zeroes in $v$ and $\ell$ the number of ones.  Then the summand corresponding to $v$ is
 \[
 \begin{minipage}{.9in}
\labellist
\small
\pinlabel $\b_v^L$ at 18 72
\pinlabel $\omega(\b_v^L)$ at 18 17
\pinlabel $K_{\ell+1}$ at 30 45
\pinlabel $k$ at 5 -3 
\pinlabel $\ell$ at  21 -3
\pinlabel $1$ at  38 -3
\pinlabel $1$ at  55 -3
\endlabellist
\begin{center}\includegraphics[scale=1]{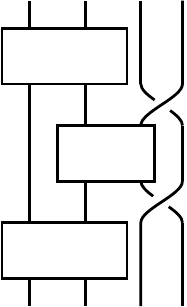}\end{center} 
\end{minipage}
 \simeq 
 \begin{minipage}{.9in}
\labellist
\small
\pinlabel $\b_v^L$ at 18 97
\pinlabel $\omega(\b_v^L)$ at 18 17
\pinlabel $K_{\ell+1}$ at 30 68
\pinlabel $k$ at 5 -3 
\pinlabel $\ell$ at  21 -3
\pinlabel $1$ at  38 -3
\pinlabel $1$ at  55 -3
\endlabellist
\begin{center}\includegraphics[scale=1]{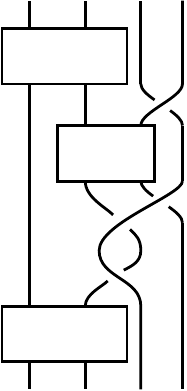}\end{center} 
\end{minipage}
\simeq
\left(
t^{\frac{1}{2}(-\ell-1)}\begin{minipage}{.9in}
\labellist
\small
\pinlabel $\b_v^L$ at 19 76
\pinlabel $\omega(\b_v^L)$ at 19 17
\pinlabel $K_{\ell+2}$ at 38 48
\pinlabel $k$ at 5 -3 
\pinlabel $\ell$ at  21 -3
\pinlabel $1$ at  38 -3
\pinlabel $1$ at  55 -3
\endlabellist
\begin{center}\includegraphics[scale=1]{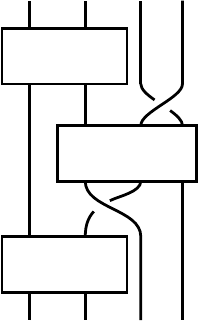}\end{center} 
\end{minipage}
\rightarrow 
q\begin{minipage}{.9in}
\labellist
\small
\pinlabel $\b_v^L$ at 19 97
\pinlabel $\omega(\b_v^L)$ at 18 17
\pinlabel $K_{\ell+1}$ at 30 68
\pinlabel $k$ at 5 -3 
\pinlabel $\ell$ at  21 -3
\pinlabel $1$ at  38 -3
\pinlabel $1$ at  55 -3
\endlabellist
\begin{center}\includegraphics[scale=1]{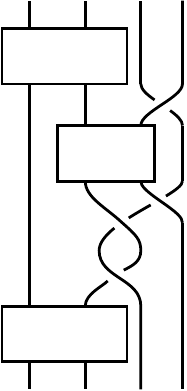}\end{center} 
\end{minipage}
\right).
 \]
The first equality (or rather, homotopy equivalence) is a simple isotopy, pulling one strand past the cable of $l$ strands. In the second equality we have used \eqref{eqn-Kn-cone}.  Applying Lemma \ref{lemma-RouquierAbsorbing} to the first complex on the right,  $K_{\ell+2}$ absorbs the $\ell$ negative crossings (in the cabled crossing below $K_{\ell+2}$) and 1 positive crossing (above), gaining an additional grading shift of $t^{\frac{1}{2}-\frac{1}{2}\ell}$.  After absorbing these crossings and applying an isotopy to the right-most complex, we obtain
\begin{equation}\label{eqn:yjminduction2}
F(\s_n)(E_{v \cdot 1} \sqcup \1_1) F(\sigma_n) \; \simeq \;  \left(
t^{-\ell}\begin{minipage}{.9in}
\labellist
\small
\pinlabel $\b_v^L$ at 19 73
\pinlabel $\omega(\b_v^L)$ at 19 15
\pinlabel $K_{\ell+2}$ at 38 46
\endlabellist
\begin{center}\includegraphics[scale=1]{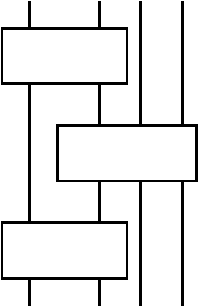}\end{center} 
\end{minipage}
\rightarrow 
q\begin{minipage}{.9in}
\labellist
\small
\pinlabel $\b_v^L$ at 19 73
\pinlabel $\omega(\b_v^L)$ at 23 15
\pinlabel $K_{\ell+1}$ at 45 44
\endlabellist
\begin{center}\includegraphics[scale=1]{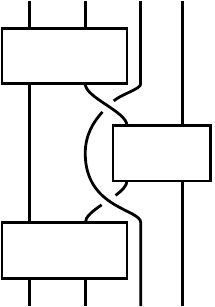}\end{center} 
\end{minipage}
\right).
 \end{equation}
The first term is just $E_{v\cdot 1 \cdot 1}$, and the second is $E_{v\cdot 0\cdot 1}$. Applying this simplification to each term of the right-hand side of (\ref{eq-yjmInduction}) completes the inductive step. It remains to verify that the grading shifts and partial order on the convolution are as claimed.

The grading shifts are determined recursively by $G_{v\cdot 1}=t^{-|v|} G_v$ and $G_{v\cdot 0}=qG_v$, from which the formula $G_v = t^{-\binom{|v|}{2}} q^{n-1-|v|}$ follows easily.

Suppose that, for two sequences $v,w \in \{0,1\}^{n-1}$, the differential from the summand $E_{v \cdot 1}$ to the summand $E_{w \cdot 1}$ is zero in the twisted differential for $F(y_n)$.  Then replacing $F(\s_n) (E_{v \cdot 1} \sqcup \1_1) F(\s_n)$ and $F(\s_n) (E_{w \cdot 1} \sqcup \1_1) F(\s_n)$ by the equivalent complexes in the right-hand side of \eqref{eqn:yjminduction2} does not introduce any differential between any of the corresponding terms in the twisted differential for $F(y_{n+1})$.  Moreover, there is no differential from $E_{v \cdot 0 \cdot 1}$ to $E_{v \cdot 1 \cdot 1}$.
Therefore, induction implies that the twisted differential respects the antilexicographic order. \end{proof}

\subsection{The full twist}

It is fairly easy to bootstrap this convolution description of the Young-Jucys-Murphy elements into a convolution description of the full twist.

\begin{defn} \label{defn-D} Let $v \in \{0,1\}^{n}$ be a shuffle, with $k$ zeroes and $\ell$ ones, so that $k + \ell = n$. Then let $D_v \in \CC^b(\SBim)$ be defined as follows:
\begin{equation} \label{eqn-D-defn} D_v := F(\b_{v})(\FT_k \sqcup K_{\ell})F(\w(\b_{v})). \end{equation} Note that the shuffle braids involved here are right-handed. \end{defn}

\begin{ex} Here is $D_v$ for $v=(10101101)$. \[
D_{10101101} =
\begin{minipage}{1.1in}
\labellist
\small
\pinlabel $\FT_3$ at 18 40
\pinlabel $K_5$ at 65 40
\endlabellist
\begin{center}\includegraphics[scale=1]{diagrams/FTsummand}\end{center} 
\end{minipage}
\]
In this example, the last index in $v$ is a one, so the rightmost strand in $D_{v}$ does not cross the strands coming from the full twist. \end{ex}

\begin{theorem}\label{thm-FTreduction}
Let $\FT_n$ denote the full right-handed twist on $n$-strands.  We have
\begin{equation}\label{eq-FTreduction}
\FT_n \simeq \bigoplus_{v\in\{0,1\}^{n-1}}q^{k} D_{v\cdot 1}
\end{equation}
with twisted differential, respecting the antilexicographic order on $\{0,1\}^{n-1}$.
Here $k$ is the number of zeroes in $v$.
\end{theorem}
We have omitted the functor $F$ from the notation, identifying a braid with its Rouquier complex.  We employ this abuse of notation frequently henceforth.
\begin{proof}
Note that the full twist braid factors as $\FT_{n+1}=\FT_{n}y_{n+1}$, where as usual $\FT_n$ is viewed as an element inside $\Br_{n+1}$ via the inclusion $\Br_n \to \Br_{n+1}$.  Actually it will be more useful to write
\[
\FT_{n+1} = \HT_{n} y_{n+1} \HT_{n}.
\]
Proposition \ref{prop-yjmSimp} gives an expression of the Jucys-Murphy complex $F(y_{n+1})$.  Tensoring on the left and right with $\HT_{n}=\HT_n\sqcup \one_1$ gives
\begin{eqnarray*}
\FT_{n+1} & \simeq &\bigoplus_v G_v (\HT_{n} \sqcup \one_1) (\b_v^L\sqcup \one_1)(\one_k\sqcup K_{\ell+1}) (\omega(\b_v^L)\sqcup \one_1)(\HT_{n}\sqcup \one_1)\\
 & \simeq &\bigoplus_v G_v ((\HT_{n}\b_v^L)  \sqcup \one_1) (\one_k\sqcup K_{\ell+1}) ((\omega(\b_v^L)\HT_{n})\sqcup \one_1)\\
  & \simeq &\bigoplus_v G_v (\b_{r(v)}  \sqcup \one_1) (\HT_k\sqcup \HT_\ell\sqcup \one_1)(\one_k\sqcup K_{\ell+1}) (\HT_k\sqcup \HT_\ell\sqcup \one_1) ((\omega(\b_{r(v)}) \one_1)\\
   & \simeq &\bigoplus_v G_v t^{\binom{\ell}{2}} (\b_{r(v)}  \sqcup \one_1) (\FT_k\sqcup K_{\ell+1}) ((\omega(\b_{r(v)}) \one_1).
\end{eqnarray*}
Here, $G_v=t^{-\binom{\ell}{2}}q^k$ is the shift determined by Proposition \ref{prop-yjmSimp}.  The second equivalence is simply given by reassociating.  The third equivalence holds by Proposition \ref{prop-HTshuffle}, which describes how $\HT_n$ interacts with shuffle braids.  The last equivalence holds since $K_{\ell+1}$ absorbs the two copies of $\HT_\ell$ (each of which has $\binom{\ell}{2}$ crossings), and the two copies of $\HT_k$ contribute a factor of $\FT_k$.  The grading shift on each summand is $G_v t^{\binom{\ell}{2}} = q^k$, as claimed.  This completes the proof.
\end{proof}

\begin{remark} One can construct convolution descriptions of other torus links in much the same way. At the moment, we have done this ad hoc for small torus links, and have neglected to
write it down here for reasons of space. It would be interesting to find a combinatorial framework (analogous to shuffles) in order to treat the general torus link in a more holistic
fashion. The results of the next chapter, including the parity miracle which makes these convolution descriptions useful, can also be adapted to our small examples in a straightforward way.
The fruits of this labor are presented in the appendix. \end{remark}

%% file: FTHHHconvolve3.tex
\section{Resolving the Hochschild homology of the full twist}
\label{sec-mainResult}

In this section we introduce Hochschild cohomology $\HH$, and we compute $\HH(\FT_n)$ for all $n\geq 1$. Our strategy is recursive. The main result of the previous section expresses
$\FT_n$ as a filtered complex whose subquotients are of the form $D_v$. In this section we show that $\HH(D_v)$ has a filtration in terms of other $\HH(D_w)$, with smaller $w$. However,
we will find it convenient to work instead with related complexes $C_v$, to be defined in section \S\ref{subsec-ourcomplexes}.

\subsection{Hochschild cohomology} \label{subsec-HHH}

The zeroth Hochschild cohomology functor $\HH^0$ is the functor which takes a graded $R$-bimodule $M$ to the graded vector space $\oplus_{m \in \ZM} \Hom(R,M(m))$ of bimodule maps of all
degrees. Its higher derived functors $\HH^k$ are packaged together in a single functor $\HH = \oplus_{k \ge 0} \HH^k \co \SBim_n \to \QM\vect^{\ZM \times \ZM}$, where this latter category
is the category of bigraded vector spaces. The two gradings are the internal grading of the bimodule (the $m$ in the direct sum above), and the Hochschild cohomological grading $k$, which
we call the \emph{Hochschild grading}.

Extending to complexes gives a functor $\HH \co \CC^b(\SBim_n) \to \CC^b(\QM\vect^{\ZM\times \ZM})$.  Given a complex $C$, $\HH^0(C)$ is the complex $\RHom(R,C)$ used to compute maps of all internal and homological degrees from the complex $R$ (concentrated in a single homological degree) to $C$. The homology of $\HH(C)$ and $\HH^0(C)$ are denoted by $\HHH(C)$ and $\HHH^0(C)$.  When we wish to emphasize the index $n$ (which is not a degree, but the number of strands) we will write $\HH(R_n;C)$, $\HH^0(R_n;C)$, and so on.

Note that $\HH(C)$ is triply graded. We will denote shifts in the tridegree by $Q^iA^jT^k \HH(C)$, where $Q$ is the usual degree, $A$ is Hochschild degree, and $T$ is homological degree.
In previous sections we found it useful to introduce the variables (or grading shifts) $t=T^2Q^{-2}$ and $q=Q^2$.  In the sequel it will prove convenient to package the Hochschild and $q$-degrees together by introducing $a=AQ^{-2}$.  We write $\PC_C(q,a,t)$ for the Poincar\'e series of $\HHH(C)$.  

The experienced reader may wish to orient himself or herself by observing that in these conventions, we have
\[
\PC_{R_1} = \frac{1+Q^{-2}A}{1-Q^2}=\frac{1+a}{1-q}
\]
which is the Hochschild cohomology of
the ring $\QM[x]$ as a bimodule over itself. The reader should think of $(1-q)\inv$ as the Poincar\'e series of $\QM[x]$ itself, and $(1 + a)$ as the Poincar\'e series of the
exterior algebra in one variable. Similarly, $\PC_{R_n} = (1+a)^n (1-q)^{-n}$.


\begin{ex} We have $\PC_{B_s} = (1-Q^2)^{-n} (1+Q^{-2}A)^{n-1} (Q+Q^{-3} A)$. Here is a brief conceptual explanation. Consider the Koszul complex which resolves $R$ by free
$R$-bimodules. Applying $\Hom$ to $R$, the differentials all become zero, yielding $(1+Q^{-2}A)^n$ times the Poincar\'e series of $R$, as for $\PC_R$ above. Applying $\Hom$ to $B_s$
instead, one of the differentials in the Kozsul complex is non-zero, becoming the middle differential in Example \ref{ex-K2}, except dualized. Thus this differential yields a factor of
$(Q+Q^{-3}A)$ instead. We will not use this computation. \end{ex}

\begin{definition}\label{def-simEquiv1} We say that two complexes $A ,B \in \CC^b(\SBim_n)$ are \emph{$\HH$-equivalent} if $\HH(A)\simeq \HH(B)$ as complexes of triply graded vector
spaces. In this case we will write $A\sim B$. \end{definition}

The basic property of Hochschild cohomology which motivates its relationship with braid closures is that $\HH(CD) \cong \HH(DC)$ whenever these tensor products make sense (e.g. if
$C,D \in \CC(\SBim_n)$ are simultaneously bounded above or below). Thus, any complex $C$ is $\HH$-equivalent to $F(\b) C F(\b\inv)$ for any braid $\b$.

Note that $\HH^i$ can actually be viewed as a map from $R$-bimodules to the subcategory of $R$-bimodules for which the left and right actions agree, which can be identified with
$R$-modules. However, the isomorphism $\HH(CD) \cong \HH(DC)$ is not an isomorphism of (complexes of) $R$-bimodules, only of their underlying vector spaces. Nonetheless, there is
still an action of $R_n$ on any Hochschild complex $\HH^i(R_n;C)$.

\begin{remark} The isomorphism $\HH(R_n;CD) \cong \HH(R_n;DC)$ of complexes of vector spaces does actually lift to an isomorphism of modules over the invariant subring $R^{S_n}$. \end{remark}

\subsection{The Markov move for Jones-Wenzl projectors} \label{subsec-markov}

The Markov move states that the closure of a braid $\b$ on $n-1$ strands is isotopic (as a link) to the closure of the braid $\s^\pm_{n-1} (\b \sqcup \1_1)$ on $n$ strands. To prove that
$\HHH$ is a link invariant, Khovanov \cite{Kh07} proved a result comparing $\HH(R_{n-1};\b)$ and $\HH(R_n;\s^\pm_{n-1} (\b \sqcup \1_1))$. In this paper, we will need a similar result, comparing
$\HH(R_{n-1};K_{n-1})$ and $\HH(R_n;K_n)$.


\begin{proposition}\label{prop-HHprops} Suppose that $2 \le n$. Let $C \in \CC^b(\SBim_{n-1})$ be viewed as a complex in $\CC^b(\SBim_n)$ via the usual inclusion functor. We have
\begin{equation} \label{eqn-HH-Kn} \HH(R_n;CK_n) \simeq t^{n-1} \HH(R_{n-1};CK_{n-1}) \oplus a\HH(R_{n-1};CK_{n-1}). \end{equation}  This can also be described as
\[
\HH(R_n;CK_n)\simeq t^{n-1} \HH(R_{n-1};CK_{n-1})\otimes \Lambda[\xi_n]
\]
where $\deg(\xi_n)= t^{1-n}a$. Hence $\PC_{CK_n} = (t^{n-1}+a)\PC_{CK_{n-1}}$. 
\end{proposition}

\begin{proof}  We use results in \cite{Hoga-pp}.  Let $\CC_n$ denote the bounded derived category of graded $(R_n,R_n)$-bimodules, where $R_n=\Q[x_1,\ldots,x_n]$ as usual.   Let $\DC_n=\K^b(\CC_n)$ denote the homotopy category of $\CC_n$.  Note that $\SBim_n$ includes as a full subcategory of $\CC_n$, and $\K^b(\SBim_n)$ includes as a full subcategory of $\DC_n$.  In case $n=0$, $\DC_0$ is equivalent to the category $\Q^{\Z\times \Z\times \Z}$ of triply graded vector spaces.

There is a \emph{partial Hochschild cohomology functor} $T_n:\DC_n\rightarrow \DC_{n-1}$, such that $\HHH=T_1\circ \cdots \circ T_n$.  These can be defined as the right adjoints to the standard inclusions $I_n:\DC_{n-1}\rightarrow \DC_n$.  We usually abuse notation, and write $C$ when we mean $I(C)$. For each $C\in \DC_{n-1}$ and each $D\in \DC_n$ we have
\[
T(CD)\cong C T(D) \ \ \ \ \ \ \ \ \text{ and } \ \ \ \ \ \ \ \ \ T(DC)\cong  T(D)C.
\]
Recall from the proof of Proposition \ref{prop:K} that $K_n\simeq K_{n-1}Q_n$, so that $T_n(K_n)\simeq K_{n-1}T_n(Q_n)$.  It was proven in \S 4 of \cite{Hoga-pp} that $T_n(Q_n)\simeq t^{n-1} P_{n-1}\otimes \Lambda[\xi_n]$.  Since $K_{n-1}P_{n-1}\simeq K_{n-1}$, we conclude that
\begin{equation}\label{eq:T(K)}
T_n(K_n)\simeq  K_{n-1}T(Q_n) \simeq t^{n-1} K_{n-1}\otimes \Lambda[\xi_n].
\end{equation}
From this, the Lemma follows from the observation that $\HHH=T_1\circ \cdots \circ T_n$.
\end{proof}

In our reduction algorithm to come, we need a relative version of the above proposition.
	
\begin{cor}\label{cor:HHHKunreduced} Suppose that $2 \le \ell \le n$, with $n = k + \ell$. Let $C \in \CC^b(\SBim_{n-1})$ be viewed as a complex in $\CC^b(\SBim_n)$ via the usual inclusion functor. We have
\begin{equation} \label{eqn-HH-Kl} \HH(R_n;C(\1_k \sqcup K_\ell)) \simeq t^{\ell-1} \HH(R_{n-1};C(\1_k \sqcup K_{\ell-1}))\otimes \Lambda[\xi_\ell]. \end{equation}
Hence $\PC_{C(\1_k \sqcup K_\ell)} = (t^{\ell-1}+a)\PC_{C(\1_k \sqcup K_{\ell-1})} $.
\end{cor}
\begin{proof}
We picture the partial trace $T_n:\DC_n\rightarrow \DC_{n-1}$ graphically as identifying the top right and bottom right strands.  The statement of the Lemma then becomes
\begin{equation}\label{eq:partialTrace}
\begin{minipage}{1.1in}
\labellist
\small
\pinlabel $K_\ell$ at 32 25
\pinlabel $C$ at 19 52
\pinlabel $k$ at 3 5
\pinlabel $\ell-1$ at 43 1
\pinlabel 1 at 59 25
\endlabellist
\begin{center}\includegraphics[scale=1]{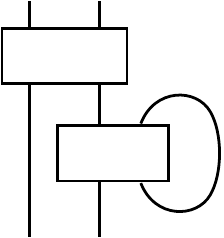}\end{center} 
\end{minipage}
\ \ \ \simeq \ \ \ t^{\ell-1}\begin{minipage}{.9in}
\labellist
\small
\pinlabel $K_{\ell-1}$ at 29 20
\pinlabel $C$ at 19 48
\pinlabel $k$ at 3 5
\pinlabel $\ell-1$ at 45 5
\endlabellist
\begin{center}\includegraphics[scale=1]{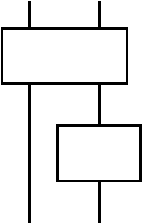}\end{center} 
\end{minipage}\otimes \Lambda[\xi_\ell]
\end{equation}
whose proof is immediate given (\ref{eq:T(K)}).
\end{proof}

This corollary allows one to slowly shrink copies of $K_\ell$ that appear, reducing the number of strands in the process. However, the case $\ell = 1$ is separate; diagrammatically, this
corresponds to taking the complex $C$ and adding a circle, since $K_1$ is the identity.

\begin{prop} Let $C \in \CC^b(\SBim_{n-1})$ be viewed as a complex in $\CC^b(\SBim_n)$ via the usual inclusion functor. We have
\begin{equation}\label{eqn-HH-K1}
\HH(R_n;C \sqcup K_1) \simeq \HH(R_{n-1};C) \ot \QM[x_n] \ot \Lambda[\xi_1]
\end{equation}
where $\deg(\xi_1) = Q^{-2} A = a$. Hence $\PC_{C \sqcup K_1} = (1-q)^{-1}(1+a)\PC_{C}$. \end{prop}

\begin{proof} In general, $\HHH(A \sqcup B) \cong \HHH(A) \ot \HHH(B)$. So this proposition just amounts to the observation that $\HHH(\QM[x_n]) \cong \QM[x_n] \ot \Lambda[\xi_1]$.  This may be pictured as
\[
\begin{minipage}{.8in}
\labellist
\small
\pinlabel $C$ at 12 24
\pinlabel $n-1$ at 27 -1
\pinlabel 1 at 56 25
\endlabellist
\begin{center}\includegraphics[scale=1]{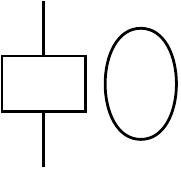}\end{center} 
\end{minipage}
\ \ \ \simeq \ \ \ t^{\ell-1}\begin{minipage}{.5in}
\labellist
\small
\pinlabel $C$ at 12 24
\pinlabel $n-1$ at 27 5
\endlabellist
\begin{center}\includegraphics[scale=1]{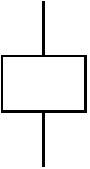}\end{center} 
\end{minipage}\otimes \Q[x_n]\otimes\Lambda[\xi_1].
\]
\end{proof}

Results like these we will also call $\HH$-equivalences.

\begin{definition}\label{def-simEquiv2} We extend Definition \ref{def-simEquiv1} above to say that two complexes $A \in \CC^b(\SBim_n)$ and $B \in \CC^b(\SBim_k)$ are
\emph{$\HH$-equivalent}, written $A \sim B$, if $\HH(A)\simeq \HH(B)$ as complexes of triply graded vector spaces. Note that $n$ and $k$ need not be equal. The complexes $A$ and $B$ are
allowed to have built-in formal Hochschild grading shifts. \end{definition}

\subsection{Reduced complexes} \label{subsec-reducedcomplexes}

Observe that the formulas for $\PC_{C(\1_k \sqcup K_\ell)}$ for $\ell \ge 2$ and $\ell = 1$ do not follow the same pattern, as the $\ell=1$ case has an extra factor of $(1-q)^{-1}$.
We will need to use both \eqref{eqn-HH-Kl} and \eqref{eqn-HH-K1} in our recursion for $\HHH(\FT_n)$, and the differences between these two formulas lead to a bookkeeping nightmare.
Instead, we will introduce the reduced complexes $\hat{K}_\ell$, which admit a streamlined formula which works for $\ell=1$ as well as for $\ell \ge 2$.

\begin{defn}\label{def:Chat}
For each complex $C\in\K(\SBim_n)$ and each element $f\in R_n$ of the ground ring, let $f\Id_C$ and $\Id_Cf$ denote the endomorphisms of $C$ given by left and right multiplication by $f$, respectively.  Set $\hat{C}:=\Cone(x_n\Id_C)$.
\end{defn}

\begin{remark}
Effectively, $\hat{C}$ is the result of killing the left action of $x_n$ on $C$.   Indeed, since $R_n$ acts freely on Soergel bimodules, standard arguments imply that $\hat{C}\simeq C/x_nC$.  Thus, one may think of $\hat{C}$ as a reduced version of $C$. The usual reduced complex is the mapping cone on $e_1\Id_C$ (equivalently, the quotient $C/e_1C)$, where $e_1 = x_1+\cdots +x_n$.  The two notions are related, but generally different.  Our sole reason for introducing $\hat{C}$ is have a functorial way of converting expressions involving $K_n$ to expressions involving $\hat{K}_n$.
\end{remark}

The relationship between $\hat{K}_n$ and $K_n$ is best understood as follows.  Let $P_n$ be the complex introduced in \cite{Hoga-pp}, whose definition is recalled in the proof of Proposition \ref{prop:K}.  Let $\END(P_n)$ denote the bigraded ring spanned by all homogeneous chain maps $Q^iT^jP_n\rightarrow P_n$ modulo homotopy.  In \cite{Hoga-pp} the second author showed that $\END(P_n)\cong \Q[u_1,u_2,\ldots,u_n]$, where the $u_k$ are variables of bidegree $\deg(u_k)= Q^{2k}T^{2-2k} = q t^{1-k}$.  In particular, $u_1$ has degree $q$, and is represented by left or right multiplication by $x_i\in R$ in the ground ring (in \cite{Hoga-pp} it is shown that all such maps are homotopic; see also the proof of Lemma \ref{lemma:K1hopping} below).

The complex $K_n$ can be interpreted as the total complex of the Koszul complex associated to the action of $u_2,\ldots,u_n$ acting on $P_n$.  Precisely: $K_n\simeq \Cone(u_2)\Cone(u_3)\cdots \Cone(u_n)$, where concatenation denotes tensor product.  This description clarifies the manner in which our definition of $K_n$ gives special treatment to the case $n=1$.  A more equitable construction would also include a factor of $\Cone(u_1)$.  By our comments above, this is precisely what $\hat{K}_n$ does: $\hat{K}_n \simeq \Cone(u_1)\Cone(u_2)\cdots \Cone(u_n)$.

\begin{lemma}\label{lemma:K1hopping}
We have $\hat{K}_n\simeq \hat{K}_1K_n$. 
\end{lemma}
\begin{proof}
It is clear that $\hat{K}_1K_n =\Cone(x_1\Id_{K_n})$, while $\hat{K}_n=\Cone(x_n\Id_{K_n})$ by definition.  Thus, it suffices to prove that $x_1\Id_{K_n}\simeq x_n\Id_{K_n}$, and the Lemma will follow by homotopy invariance of mapping cones.  It is a standard fact that there are canonical maps $R(-1)\rightarrow B_i$ and $B_i(-1)\rightarrow R$ whose composition is $\a_i:=x_i-x_{i+1}$.  Thus, $B_iK_n \simeq 0$ implies that $\a_i\Id_{K_n}$ factors through a contractible complex, hence is null-homotopic for all $1\leq i\leq n-1$.  This implies that $x_1\Id_{K_n}\simeq x_n\Id_{K_n}$, and completes the proof.
\end{proof}

\begin{remark}
Applying the functor $C\mapsto \hat{C}$ to the result of Theorem \ref{thm-FTreduction} yields an equivalence $\hat{\FT}_n \simeq \bigoplus_{v\in\{0,1\}^{n-1}}q^{k} \hat{D}_{v\cdot 1}$ with twisted differential.  Further, $\hat{D}_{v\cdot 1}$ is given by the same formula as $D_{v\cdot 1}$, except with $K_\ell$ replaced by $\hat{K}_\ell$.
\end{remark}

Our next result will later be used to show that the computation of $\HHH(\FT_n)$ reduces to a computation of $\HHH(\hat{\FT}_n)$.  In particular $\PC_{\FT_n} = \frac{1}{1-q}\PC_{\hat{\FT}_n}$.

\begin{prop}\label{prop-reducedHHH}
If $C\in \CC^b(\SBim_n)$ is such that $\HHH(\hat{C})$ is supported in even homological degrees, then so is $\HHH(C)$, and $\HHH(C)\cong \Q[x_n]\otimes \HHH(\hat{C})$.  In particular, if $\PC_C(q,a,t)$ denotes the Poincar\'e series of $\HHH(C)$, then
\[
\PC_C = \frac{1}{1-q} \PC_{\hat{C}}.
\]
\end{prop}
\begin{proof} Consider a more general situation in which $M$ is a chain complex on which some polynomial ring $\Q[x]$ acts.  Let $Z$ denote the dg algebra $\Q[x,y,\theta]$ with $d(\theta)=x-y$, $d(x)=0$, and $d(y)=0$.  Here, $\theta$ is an odd variable, hence we assume that $\theta^2=0$.  The differential ensures that $y\simeq x$.  More precisely, there is a chain map $Z\rightarrow \Q[x]$ sending $\theta\mapsto 0$, $x\mapsto x$, and $y\mapsto x$.  This map is a homotopy equivalence $Z\rightarrow \Q[x]$.  Further, the inverse map and the relevant homotopies can all be chosen to be $\Q[x]$-equivariant.

Consider the chain complex $M'\simeq Z\otimes_{\Q[x]} M$.  This is regarded as a dg $\Q[x]$-module in a slightly non-standard way, where $x$ acts by multiplication by $y$ on the first tensor factor.  The above paragraph implies that there is a homotopy equivalence $M'\simeq M$ which commutes the $\Q[x]$ actions up to homotopy. Now, $M' \cong \Q[y,\theta]\otimes M$ with twisted differential
\begin{itemize}
\item $d(y\otimes m)  = y\otimes d(m)$,
\item $d(\theta\otimes m) = y\otimes m -1\otimes x m - \theta\otimes d(m)$. 
\end{itemize}
After rearranging, we see that $M\simeq M'\simeq \Q[y]\otimes \hat{M}$ with twisted differential, where $\hat{M}$ denotes the ``reduced complex'' $\hat{M}:=\Cone(M\buildrel x\over\to M)$.  This construction is formally analogous to the fact that if $X$ is a topological space on which a group $G$ acts, then there is a space $X'$ on which $G$ acts \emph{freely}, such that $X\simeq X'$ via a $G$-equivariant homotopy equivalence.

Now we apply this construction to the case of interest.  Since $\HH$ is a linear functor which is extended to complexes, we have that $\HH$ commutes with mapping cones.  In particular $\HH(\hat{C})$ is the mapping cone of $x_n$ acting on $\HH(C)$.  The above construction then produces a twisted differential on $\Q[x_n]\otimes \HH(\hat{C})$ such that the resulting complex is homotopy equivalent to $\HH(C)$.  Thus we may regard $\HH(C)$ as a convolution of complexes $q^k\HH(\hat{C})$, indexed by $k\in \Z_{\geq 0}$.  The result now follows by the parity miracle (Proposition \ref{prop:evenConvolutions}).  Strictly speaking the parity miracle doesn't directly apply, because the indexing set is not finite.  To fix this problem we fix $r$ and consider the subcomplex of $\HH(C)$ consisting of chains with $q$-degree $r$. Each of these is a convolution with only finitely many terms, since $q^k\HH(\hat{C})$ is supported in large $q$-degrees for $k$ large.  Thus, the parity miracle can be applied separately to each $q$-degree.  The details are straightforward, so we omit them.
 \end{proof}

\begin{proposition}\label{prop:Khat}
The complexes $\hat{K}_n$ satisfy the same recursion as $K_n$.  That is:
\begin{equation} \label{eqn-Khat-cone}
\begin{minipage}{1in}
\labellist
\small
\pinlabel $\hat{K}_{n-1}$ at 27 39
\endlabellist
\includegraphics[scale=1]{diagrams/KwithYJM0}
\end{minipage}
\ \ \simeq \ \ \left( \
\begin{minipage}{1in}
\labellist
\small
\pinlabel $\hat{K}_{n}$ at 28 29
\endlabellist
\includegraphics[scale=1]{diagrams/KwithYJM1}
\end{minipage}
\longrightarrow
q \ \begin{minipage}{1in}
\labellist
\small
\pinlabel $\hat{K}_{n-1}$ at 28 29
\endlabellist
\includegraphics[scale=1]{diagrams/KwithYJM2}
\end{minipage}
\right).
\end{equation}
This holds for all $n\geq 1$, where by convention we set $\hat{K}_0=\Q\in \SBim_0$.
\end{proposition}
\begin{proof}
For $n\geq 2$ simply tensor \eqref{eqn-YJM-cone} on the left with $\hat{K}_1$ and use Lemma \ref{lemma:K1hopping}.  For the somewhat degenerate case $n=1$, the result follows from the following argument.  Note that when $n=1$, the left hand side of \eqref{eqn-Khat-cone} is $\Q[x_1]$, and the second term on the right-hand side can be identified with $x_1\Q[x_1]$.  Then the result follows from the observation that $\hat{K}_1\simeq \Q$, and $\Q[x_1]\cong \Q \oplus x_1\Q[x_1]$.
\end{proof}

We have the following streamlined version of the Markov move for the reduced Jones-Wenzl complexes.

\begin{cor} \label{cor-streamlined} Suppose that $1 \le \ell \le n$, with $n = k + \ell$. Let $C \in \CC^b(\SBim_{n-1})$ be viewed as a complex in $\CC^b(\SBim_n)$ via the usual inclusion functor. We have
\begin{equation} \label{eqn-HH-Klhat} \HH(R_n;C(\1_k \sqcup \hat{K}_\ell)) \simeq  t^{\ell-1}\HH(R_{n-1};C(\1_k \sqcup \hat{K}_{\ell-1})) \otimes \Lambda[\xi_\ell] \end{equation}
where $\Lambda$ denotes an exterior algebra, and $\deg \xi_i = t^{1-i} a$.  Hence $\PC_{C(\1_k \sqcup \hat{K}_\ell)} =(t^{\ell-1}+ a) \PC_{C(\1_k \sqcup \hat{K}_{\ell-1})} $. \end{cor}
\begin{proof}
For $\ell>1$ this follows by the same argument in the proof of Corollary \ref{cor:HHHKunreduced}.  For $\ell=1$ this follows from the fact that $\HH(\hat{K}_1)$ is the mapping cone of $x_1$ acting on $\HH(\Q[x_1];\Q[x_1])=\Q[x_1]\otimes \Lambda[\xi_1]$.  Since $x_1$ acts freely, standard arguments imply that this mapping cone is equivalent to the quotient $\Lambda[\xi_1]$, which has Poincare polynomial $1+a$.
\end{proof}

\begin{cor} \label{cor-hatKn} Let $\Lambda = \Lambda[\xi_1, \ldots, \xi_n]$ be the exterior algebra.  Let $\PC_\Lambda$ denote its Poincar\'e polynomial. Then $\PC_{\hat{K}_n} = t^{\binom{n}{2}}\PC_{\Lambda}$. \end{cor}

Note that
\begin{equation} \PC_{\hat{K}_n} = \prod_{i=1}^n (t^{i-1} + a) = \prod_{i=0}^{n-1} (t^i + a). \end{equation}

\subsection{The complexes we use} \label{subsec-ourcomplexes}

Now let us return to the convolution description of the full twist.

\begin{defn} \label{def-Cvprime} For each $v\in \{0,1\}^n$ with $k$ zeroes and $\ell$ ones, let $C_v = \Tw_v (\FT_k \sqcup \hat{K}_\ell)$. Recall that the shuffle twist $\Tw_v$ was
described in Definition \ref{def-shuffleBraids}. Again, we identify a braid with its Rouquier complex. \end{defn}

\begin{example}
If $v = (0101100)$, then we have
\[
C_v \ = \ \begin{minipage}{1.5in}
\labellist
\small
\pinlabel $\FT_4$ at 20 20
\pinlabel $\hat{K}_3$ at 70 20
\endlabellist
\begin{center}\includegraphics[scale=1]{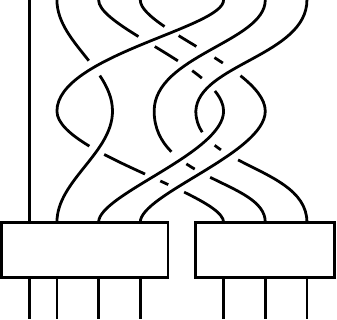}\end{center} 
\end{minipage}
\]
\end{example}

It is clear that $C_v$ is conjugate to $\hat{D}_v$, where $D_v$ is as in Definition \ref{defn-D}.  These will be the complexes we use in our inductive computation of Hochschild cohomology.


\begin{example}
We have $C_{00 \cdots 0} = \FT_n$ while $C_{10 \cdots 0} = \hat{\FT}_n$.  In general $C_{1 \cdots 10 \cdots,0}$ is given by a diagram of the form
\[
C_{1111000} \ = \ \ig{1}{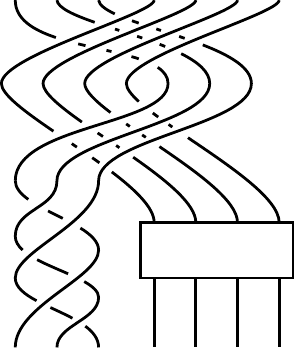}
\]
exemplified here in the case of $v=(1111000)$. The empty white box represents $\hat{K}_4$.
\end{example}

We now present a simple convolution recursion for $\HH(C_v)$.

\subsection{The key recursion} \label{subsec-newRecursion}

\begin{lemma}\label{lemma-shuffleRecursion2}
Let $v\in \{0,1\}^n$ be a sequence with $k$ zeroes and $\ell = n-k$ ones.  Then
\[
\Tw_{v\cdot 0} \ \ = \ \ 
\begin{minipage}{1.5in}
\labellist
\small
\pinlabel $k$ at 7 6
\pinlabel $1$ at 26 6
\pinlabel $\ell$ at 57 6
\pinlabel $\Tw_v$ at 34 51
\endlabellist
\begin{center}\includegraphics[scale=.7]{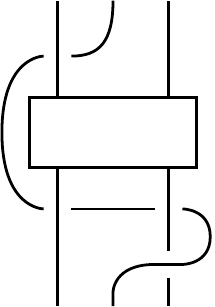}\end{center} 
\end{minipage}.
\]
\end{lemma}
\begin{proof}
This is simply the result of Proposition \ref{prop-shuffleRecursion}, followed by an isotopy.
\end{proof}

\begin{proposition}\label{prop:newRecursion}
We have $\HH(C_{v\cdot 0}) \simeq \Big(\HH(C_{1\cdot v})\rightarrow q\HH(C_{0\cdot v})\Big)$.
\end{proposition}

\begin{proof}
Let $v\in\{0,1\}^n$ be a sequence with $k$ zeroes and $\ell = n-k$ ones. Observe that $C_{v\cdot 0}$ can be written
\[
C_{v\cdot 0} \ \  = \ \ 
\begin{minipage}{1.1in}
\labellist
\small
\pinlabel $k$ at 5 6
\pinlabel $1$ at 38 6
\pinlabel $\ell$ at 72 6
\pinlabel $\FT_{k+1}$ at 24 23
\pinlabel $K_{\ell}$ at 67 44
\pinlabel $\Tw_v$ at 41 100
\endlabellist
\begin{center}\includegraphics[scale=.8]{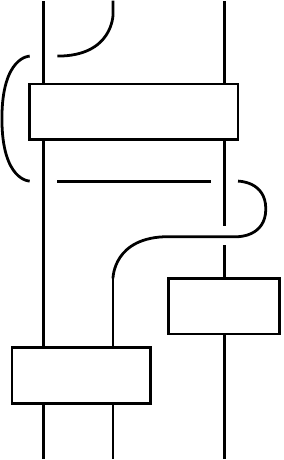}\end{center}
\vskip2pt
\end{minipage},
\]
where we have used Lemma \ref{lemma-shuffleRecursion2} to rewrite $\Tw_{v\cdot 0}$.  Now, the exact triangle for $K_n$ \eqref{eqn-YJM-cone} says that this complex is homotopy equivalent to a convolution of the form
\[
C_{v\cdot 0} \ \ \simeq \ \ \left(
\begin{minipage}{1.1in}
\labellist
\small
\pinlabel $k$ at 5 6
\pinlabel $1$ at 38 6
\pinlabel $\ell$ at 72 6
\pinlabel $\FT_{k+1}$ at 24 23
\pinlabel $K_{\ell+1}$ at 52 48
\pinlabel $\Tw_v$ at 41 93
\endlabellist
\begin{center}\includegraphics[scale=.8]{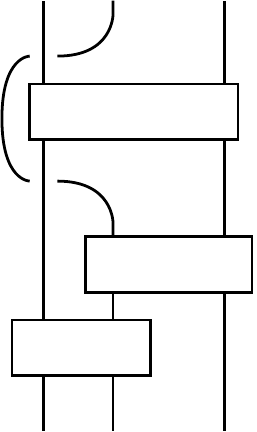}\end{center} 
\end{minipage}
\ \ \ \longrightarrow \ \ \ 
q\begin{minipage}{1.1in}
\labellist
\small
\pinlabel $k$ at 5 6
\pinlabel $1$ at 38 6
\pinlabel $\ell$ at 73 6
\pinlabel $\FT_{k+1}$ at 24 24
\pinlabel $K_{\ell}$ at 67 24
\pinlabel $\Tw_v$ at 41 69
\endlabellist
\begin{center}\includegraphics[scale=.8]{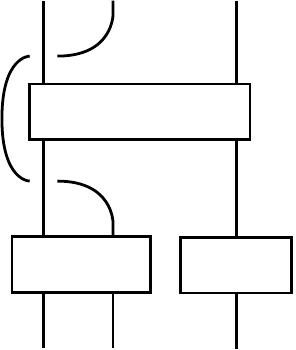}\end{center} 
\end{minipage}
\right).
\]
To prove the Proposition, we must show that the term on the left is $\HH$-equivalent to $C_{1\cdot v}$, and the term on the right is $\HH$-equivalent to $C_{0\cdot v}$.    For the term on the right, simply slide the left-handed crossing (rather, the cabled crossing between 1 strand and $k$ parallel strands) from the top to the bottom, through the full twist, where it meets and annihilates the right-handed crossing.  The resulting complex is $C_{0\cdot v}$.  For the term on the left, we have the following sequence of simplifications:
\[
\begin{minipage}{1.1in}
\labellist
\small
\pinlabel $k$ at 5 6
\pinlabel $1$ at 38 6
\pinlabel $\ell$ at 72 6
\pinlabel $\FT_{k+1}$ at 24 23
\pinlabel $K_{\ell+1}$ at 52 48
\pinlabel $\Tw_v$ at 41 93
\endlabellist
\begin{center}\includegraphics[scale=.8]{diagrams/newRecursion2}\end{center} 
\end{minipage}
\ \ \ = \ \ \ 
\begin{minipage}{1.1in}
\labellist
\small
\pinlabel $k$ at 5 6
\pinlabel $1$ at 46 6
\pinlabel $\ell$ at 72 6
\pinlabel $\FT_{k}$ at 17 44
\pinlabel $K_{\ell+1}$ at 53 68
\pinlabel $\Tw_v$ at 41 112
\endlabellist
\begin{center}\includegraphics[scale=.8]{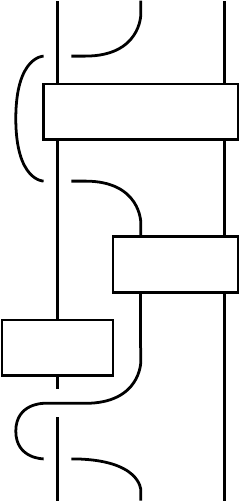}\end{center} 
\end{minipage}
\ \ \ \sim \ \ \ 
\begin{minipage}{1.1in}
\labellist
\small
\pinlabel $k$ at 5 6
\pinlabel $1$ at 47 6
\pinlabel $\ell$ at 72 6
\pinlabel $\FT_{k}$ at 17 20
\pinlabel $K_{\ell+1}$ at 52 44
\pinlabel $\Tw_v$ at 41 88
\endlabellist
\begin{center}\includegraphics[scale=.8]{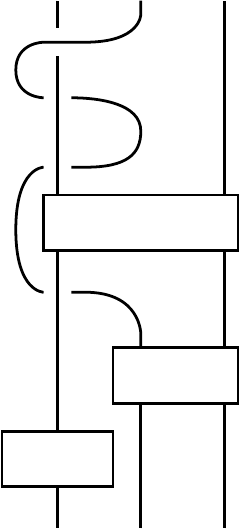}\end{center} 
\end{minipage}
\ \ \ \simeq \ \ \ 
\begin{minipage}{1.1in}
\labellist
\small
\pinlabel $k$ at 5 6
\pinlabel $1$ at 47 6
\pinlabel $\ell$ at 72 6
\pinlabel $\FT_{k}$ at 17 20
\pinlabel $K_{\ell+1}$ at 52 44
\pinlabel $\Tw_v$ at 41 88
\endlabellist
\begin{center}\includegraphics[scale=.8]{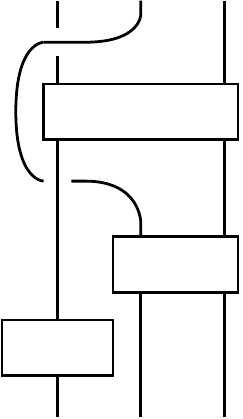}\end{center} 
\end{minipage}.
\]
The first equivalence is obtained by rewriting the full-twist as $\FT_{k+1} = \FT_k J_{k+1}$, where $J_{k+1}$ is the Jucys-Murphy braid.  The second is an $\HH$-equivalence which slides the Jucys-Murphy braid from the bottom to the top.  The final equivalence is an obvious isotopy.  The braid on the top of the resulting complex is $\Tw_{1\cdot v}$ (similar to the statement of Proposition \ref{prop-shuffleRecursion} for $\Tw_{v\cdot 0}$).  The resulting complex is therefore $C_{1\cdot v}$.  Each of the above equivalences corresponds to an honest homotopy equivalence after applying the functor $\HH(-)$. This completes the proof.  
\end{proof}

\subsection{Our main result} \label{subsec:mainResult}

In this section we prove our main theorem, which gives a recursion formula for the Poincar\'e series for $\HHH(C_v)$.
\begin{prop}\label{prop:fv}
There is a unique family of rational functions $f_v(q,a,t)$, indexed by binary sequences $v\in \{0,1\}^n$ with $n\in \Z_{\geq 0}$, satisfying
\begin{subequations}
\begin{equation} \label{eqn:recursion1} f_{v\cdot 1}(q,a,t) = (t^{|v|}+a)f_v(q,a,t)\end{equation}
\begin{equation} \label{eqn:recursion2} f_{v\cdot 0}(q,a,t) = q f_{0\cdot v} +  f_{1\cdot v}(q,a,t)\end{equation}
\end{subequations}
together with $f_{\emptyset} = 1$.
\end{prop}
\begin{proof}
Let us first prove uniqueness. Note that rule (3) applied to the sequence $v=(00\cdots 0)$ yields $f_{00\cdots 0} = qf_{00\cdots0}+f_{10\cdots 0}$.  In other words:
\begin{equation}\label{eqn:allZeroes}
f_{00\cdots0} = \frac{1}{1-q} f_{10\cdots 0},
\end{equation}
which we think of as our replacement for rule (3) when applied to the zero sequence.

Now, define a partial order on the set of binary sequences as follows: given $v\in \{0,1\}^n$ and $w\in \{0,1\}^m$, write $v< w$ if one of the following conditions is met:
\begin{itemize}
\item $n<m$
\item $n=m$ and $v$ has fewer zeroes than $w$.
\item $n=m$, $v$ and $w$ have the same number of zeroes, and number of inversions in the shuffle permuation $\pi_v$ is less than the number of inversions in $\pi_w$.
\end{itemize}
Then $<$ defines a partial order on the set of binary sequences.  The number of inversions in $\pi_v$ is the number of pairs $(i<j)$ where $v_i=1$ and $v_j=0$.  Thus, $0\cdot v \leq v\cdot 0$ with equality if and only if $v=(00\cdots 0)$.  Clearly the empty sequence is the unique minimum with respect to this partial order, and rules (2) and (3) express any $f_v$ in terms of $f_w$ with $w<v$.  This proves uniqueness.

For existence, we need to prove consistency of the rules (1), (2), and (3).  However, this is clear since for each $v$, exactly one of the rules (1), (2), (3), can be applied.
\end{proof}

We now have our main theorem:

\begin{thm}\label{thm:parityAndRecursion}
The Poincar\'e series for $\HHH(C_v)$ equals the rational function $f_v(q,a,t)$ from Proposition \ref{prop:fv}.  In particular, $\HHH(C_v)$ is supported in even homological degrees.  As a special case we have that $\HHH(\FT_n) = f_{00\cdots 0}(q,a,t)$ is the Poincar\'e series of the triply graded homology of the $(n,n)$ torus link, up to an overall shift.
\end{thm}
\begin{proof}
Let $\PC_v = \PC_{C_v}$ denote the Poincar\'e series for $\HHH(C_v)$.  We will show that $\PC_v=f_v$ by induction on $v$, using the partial order on the set of binary sequences introduced in the proof of Proposition \ref{prop:fv}.

In case $v=\emptyset$, we have $C_\emptyset = \QM$, which satisfies $\HHH(R_0; \QM) = \HHH(\QM;\QM) = \QM$.  Thus $\PC_{\emptyset} = 1 = f_{\emptyset}$.

Now, fix $v\in \{0,1\}^n$, and assume by induction that $\PC_w = f_w$ for all $w<v$.  If $v=(00\cdots 0)$, then $\PC_{10\cdots 0}=f_{10\cdots 0}$ by induction.  Further, $C_{10\cdots 0} = \hat{C}_{00\cdots 0}$, so Proposition \ref{prop-reducedHHH} says that $\PC_{00\cdots 0}=\frac{1}{1-q} \PC_{10\cdots 0}$, which equals $f_{00\cdots 0}$ by \eqref{eqn:allZeroes}.  This takes care of the case where $v$ is the zero sequence.  Thus, we assume below that $v$ is nonzero.

There are two cases: either $v=w\cdot 1$ or $w\cdot 0$ for some $w$.  In the first case, then Corollary \ref{cor-streamlined} says that $\PC_{w\cdot 1}=(t^{|v|}+a)\PC_{w}$, hence $\PC_{w\cdot 1}=f_{w\cdot 1}$ by Equation \eqref{eqn:recursion1} and induction.  Thus, we may assume that $v=w\cdot 0$.  Since $v$ is nonzero, we have $1\cdot w<w\cdot 0$ (fewer zeroes) and $0\cdot w< w\cdot 0$ (fewer inversions).  Thus, by induction, we have $\PC_{1\cdot w}=f_{1\cdot w}$ and $\PC_{0\cdot w}=f_{0\cdot w}$.  Also by induction, we may assume that $\HHH(C_{1\cdot w})$ and $\HHH(C_{0\cdot w})$ are supported in even homological degrees.  Thus, the terms of the convolution Proposition \ref{prop:newRecursion} have the same parity after taking $\HHH$. Proposition \ref{prop:evenConvolutions} implies that $\HHH(C_{w\cdot 0})$ splits as a direct sum
\[
\HHH(C_{w\cdot 0}) \cong q\HHH(C_{0\cdot w}) \oplus \HHH(C_{1\cdot w}).
\] 
Taking Poincar\'e series, we see that $\PC_{w\cdot 0}=q\PC_{0\cdot w}+\PC_{1\cdot w}$.  It follows that $\PC_v=f_v$, by induction and the uniqueness statement of Proposition \ref{prop:fv}. 
\end{proof}

%% file: FTHHHnumerology.tex
\section{Numerological considerations} \label{sec-numerology}

Below, we give an alternate recursive formula for the power series $f_v(q,a,t)$ for sequences $v \in \{0,1\}^n$.  We then give a closed formula for $a=0$ specialization $f_{00\cdots 0}(q,0,t)$.

\subsection{An alternate recursive formula}

The recursion described in this section was actually discovered before the recursion that appears in Proposition \ref{prop:fv}.  We originally proved our main result (Theorem \ref{thm:parityAndRecursion}) using this recursion, and later found a much more elegant route which now appears in our \S \ref{subsec:mainResult}.  Nonetheless this alternate recursion is quite useful, and will lead us to a closed formula for the $a$-degree zero part of $\HHH(\FT_n)$ in \S \ref{subsec:closedForm}.   We introduce some combinatorial notions which will be relevant below.

\begin{definition}\label{def-insertingW}  Fix a sequence $v \in \{0,1\}^n$.  We will call a pair of sequences $(v,w)$ \emph{compatible} if $w\in\{0,1\}^k$, where $k$ is the number of zeroes in $v$.  If $(v,w)$ is a compatible pair, we define a sequence $v\circ w$ by ``inserting $w$ into the zeroes of $v$.'' That is, let $(v \circ w)_i = 1$ if either $v_i = 1$, or if $v_i$ is the $j$-th zero in $v$ and $w_j = 1$.   We let $I_{v,w}\subset \{1,\ldots,n\}$ denote the subset of indices such that $v_i=0$ but $(v\circ w)_i=1$.   We say that $i$ is a \emph{one in $v$} if $v_i=1$, and $i$ is a \emph{one in $w$ (relative to $v$)} if $v_i=0$ but $(v\circ w)_i=1$.  When $v$ is understood, we omit the phrase \emph{relative to $v$}.
\end{definition}

\begin{example}\label{ex-inserting}
Let $v=(1101001)$ and $w=(001)$.  In this case we have $v\circ w = (11\underline{0}1\underline{0}\underline{1}1)$, where the underlined terms indicate where we have inserted $w$ into $v$.
\end{example}

\begin{defn} \label{defn-gradings} Fix $v \in \{0,1\}^n$ and $w \in \{0,1\}^k$ as above.  For each index $i$, let $\ell(i)$ denote the number of ones of $v$ strictly to the left of $i$, and let $m(i)$ denote the number of ones in $w$ strictly to the right of $i$.  Let $P_{v,w,i}$ denote $(t^{\ell(i)+m(i)}+a)$ if $v_i=1$ and $P_{v,w,i}=1$ otherwise.
\end{defn}

\begin{lemma}\label{lemma:Pvw}
The following relations hold:
\begin{enumerate}
\item $\PC_{v\cdot 1, w} = (t^{|v|}+a)\PC_{v,w}$
\item $\PC_{v\cdot 0, w\cdot 0}=\PC_{v,w}$
\item $\PC_{1\cdot v,w}=(t^{|w|}+a)\PC_{v\cdot 0, w\cdot 1}$.
\item $\PC_{0\cdot v, 0 \cdot  w}=\PC_{v,w}$
\item $\PC_{0\cdot v, 1\cdot w}=\PC_{v,w}$
\end{enumerate}
\end{lemma}
\begin{proof}
These are easily verified directly from the definition.
\end{proof}

\begin{prop}\label{prop:altRecursion} The functions $f_v(q,a,t)$ defined in Proposition \ref{prop:fv} can also be defined by the recursion:
\begin{subequations}
\begin{equation} \label{eqn-alt1} f_{00\cdots 0}(q,a,t) = (1-q)\inv f_{10\ldots0}(q,a,t)\end{equation}
\begin{equation} \label{eqn-alt2} f_v(q,a,t) = \sum_{w\in \{0,1\}^k} \PC_{v,w}(q,a,t) q^{k-|w|} f_w(q,a,t) \ \ \ \ \ \ \ (\text{ if $v\neq 0$})\end{equation}
\end{subequations}
The base of the recursion is still $f_{\emptyset} = 1$.
\end{prop}

Note that we regard $\{0,1\}^0$ as the set containing the empty set.  Thus, \eqref{eqn-alt2} gives
\[
f_{11\cdots 1}(q,a,t) = P_{11\cdots 1,\emptyset} (q,a,t) = \prod_{i=1}^n(t^{i-1}+a)
\] as a special case.

\begin{proof} Both recursions uniquely pin down a collection of functions $f_v(q,a,t)$. Therefore, if one of these definitions satisfies the other's recursive formula, then they are equivalent definitions. The recursions are clearly equivalent when computing $f_v$ for sequences of length $n \le 1$.  Let us temporarily denote by $g_v$ the family of functions determined by \eqref{eqn-alt1} and \eqref{eqn-alt2}.   We will show that the $g_v$ also satisfy the recursion which defines $f_v$ (Equations \eqref{eqn:recursion1} and \eqref{eqn:recursion2}).

First, note that $g_{v\cdot 1} =\sum_{w\in\{0,1\}^k}\PC_{v\cdot 1, w} q^{k_1} g_w$ where $k_1$ is the number of zeroes in $w$.  By part (1) of Lemma \ref{lemma:Pvw}, we have $\PC_{v\cdot 1,w}= (t^{|v|}+a)g_w$, which implies that $g_{v\cdot 1}=(t^{|v|}+a)g_v$.  Thus, the rule \eqref{eqn:recursion1} is satisfied by $g_v$.

We now check that $g_{v\cdot 0}$ satisfies  \eqref{eqn:recursion2}.  If $v=(00\cdots 0)$, then this translates precisely to rule \eqref{eqn-alt1}, which we are assuming is valid.  Thus, we may assume that $v$ is nonzero.  We must show that $g_{v\cdot 0} = q g_{0\cdot v}+g_{1\cdot v}$.    To do this, fix $v\in \{0,1\}^n$, let $\ell=|v|$ the number of ones in $v$ and $k = n-\ell$ the number of zeroes.  Below, we let $k_1$ denote the number of zeroes in $w$, so that
\[
g_v = \sum_w \PC_{v,w} q^{k_1} g_w
\]
Let us expand $g_{v \cdot 0}(q,a,t)$ using the rule \eqref{eqn-alt2}. Because of the extra zero, there are twice as many terms in this sum as there were for $g_v$, corresponding to sequences $w \cdot 0$ and sequences $w \cdot 1$. We obtain
\[
g_{v \cdot 0} = \sum_w  q^{k_1}\bigg(q \PC_{v\cdot 0,w\cdot 0}  g_{w \cdot 0} + \PC_{v\cdot 0,w\cdot 1} g_{w\cdot 1}\bigg).
\]
For the sequence $w \cdot 0$ relative to $w$, there is an extra zero yielding an extra factor of $q$.  By Lemma \ref{lemma:Pvw}, we have $\PC_{v\cdot 0,w\cdot 0}=\PC_{v,w}$, and by rule \eqref{eqn-alt1}, we have $g_{w\cdot 1}=(t^{|w|}+a)g_{w}$.  This yields
\begin{equation} \label{eqn-foobar}
g_{v \cdot 0} = \sum_w  q^{k_1}\bigg(q \PC_{v,w}  g_{w \cdot 0} + (t^{|w|}+a)\PC_{v\cdot 0,w\cdot 1} g_{w}\bigg).
\end{equation}

Recall that we are assuming $v$ is not the zero-sequence, hence $k<n$.  Now we apply induction on $n$, assuming \eqref{eqn:recursion1} and \eqref{eqn:recursion2}  hold for $g_w$.   Applying both of these equations to the right hand side of \eqref{eqn-foobar} we obtain
\begin{equation} \label{eqn-barfoo} g_{v \cdot 0} = \sum_w q^{k_1}\bigg(q^2 \PC_{v,w} g_{0 \cdot w} + q \PC_{v,w} g_{1 \cdot w} +  (t^{|w|}+a) \PC_{v\cdot 0,w\cdot 1} g_w\bigg). \end{equation}

Similarly,
\[
g_{0 \cdot v} = \sum_w  q^{k_1}\bigg(q \PC_{0\cdot v,0\cdot w} g_{0 \cdot w} + \PC_{0\cdot v,1\cdot w} g_{1 \cdot w}\bigg)
\]
and
\[ g_{1 \cdot v} = \sum_w \PC_{1\cdot v, w} q^{k_1} g_w. \]
Lemma \ref{lemma:Pvw} says $\PC_{0\cdot v,0\cdot w}=\PC_{v,w}=\PC_{0\cdot v,1\cdot w}$ and $\PC_{1\cdot v,w}=(t^{|w|}+a)\PC_{v\cdot 0,w\cdot 1}$.  Thus,
\begin{equation} \label{eqn-barfoobar}  q g_{0 \cdot v} + g_{1 \cdot v} = \sum_w q^{k_1}\bigg(q^2 \PC_{v,w}g_{0 \cdot w} + q \PC_{v,w} g_{1 \cdot w} +(t^{|w|}+a)\PC_{v\cdot 0, w\cdot 1} g_w\bigg). \end{equation}

Comparing \eqref{eqn-barfoo} to \eqref{eqn-barfoobar}, we have shown that $g_{v\cdot 0} = q g_{0\cdot v}+g_{1\cdot v}$, as desired. \end{proof}

\subsection{The redundancy of rule \eqref{eqn-alt1}}

\begin{prop} \label{prop-rule1} Equation \eqref{eqn-alt1} follows from a verbatim application of Equation \eqref{eqn-alt2}  to the case of the zero sequence $v = (00\cdots 0)$. \end{prop}

\begin{proof} Let $v$ be the zero sequence of length $n$. When we expand $f_{00\cdots 0}$ using Equation \eqref{eqn-alt2} , we obtain a sum of $2^n$ terms, indexed by sequences $w \in \{0,1\}^n$. We claim that
the sum of the terms with $w$ ending in $1$ is actually just $f_{10\cdots 0}$; that the sum of the terms with $w$ ending in $10$ is actually just $qf_{10\cdots 0}$; the sum of the terms
with $w$ ending in $100$ is $q^2 f_{10\cdots 0}$; and so forth. Of course, the unique term where $w$ ends in $10 \cdot 0$ (length $n$) is just $q^{n-1} f_{10 \cdots 0}$, because $q^{n-1}$ is easy observed to be the coefficient of $f_{10 \cdots 0}$ in this expansion.

Given this claim, we have \[f_{00\cdots 0} = q^n f_{00\cdots 0} + (1+q+q^2 + \ldots + q^{n-1}) f_{10\cdots 0},\] which immediately implies Equation (\ref{eqn-alt1}). So it is enough to show the claim.

The claim is proven by an easy induction, using Equation \eqref{eqn:recursion2}.  For instance, 
\[
q^{n-2}f_{10\cdots 0}=q^{n-1}f_{010\cdots 0} + q^{n-2}f_{110\cdots 0}.
\]
by one application of Equation \eqref{eqn:recursion2}.  These are the two terms which end in $(10\cdots 0)$ (length $n-1$).  This proves one statement of the claim.  For the next, we apply Equation \eqref{eqn:recursion2} again and decrease the power of $q$, obtaining
\[
q^{n-3}f_{10\cdots 0} = q^{n-1}f_{0010\cdots 0} + q^{n-2}(f_{1010\cdots 0}+f_{0110\cdots 0}) + q^{n-3}f_{1110\cdots 0}.
\]
These are the 4 terms which end in $10\cdots 0$ (length $n-2$), and each appears with the correct power of $q$.  This proves the second statement of our claim.  The remaining parts of the claim follow by repeating this argument.  We leave the details to the reader. \end{proof}

\subsection{The closed form of $\HHH^0(\FT_n)$}
\label{subsec:closedForm}

One useful consequence of the alternate recursion is that it leads to a simple derivation of Theorem \ref{thm-degreezerosolve-intro}. To remind the reader, this theorem stated that
\[
f_{00\cdots 0}(q,0,t)=\sum_{\sigma}t^{a(\sigma)+b(\sigma)}q^{c(\sigma)}
\]
where the sum is over functions $\sigma:\{1,\ldots,n\}\rightarrow \Z_{\geq 0}$, and the integers $a(\sigma)$, $b(\sigma)$, $c(\sigma)$ are defined by
\begin{enumerate}
\item $a(\sigma)=\sum_{k\geq 0}\binom{|\sigma\inv(k)|}{2}$
\item $b(\sigma)$ is the number of pairs $(i,j)\in \{1,\ldots,n\}$ such that $i<j$ and $\sigma(j)=\sigma(i)+1$.
\item $c(\sigma)=\sum_{i=1}^n\sigma(i)$.
\end{enumerate}

\begin{definition}\label{def-crossings} Recall Definition \ref{def-insertingW}.  If $(v,w)$ is a compatible pair of sequences, let $c(v,w)$ denote the number of pairs of indices $i<j$ such that $i$ is a one in $v$ and $j$ is a one in $w$ (relative to $v$).
\end{definition}

\begin{example}\label{ex-crossings}
Let $v=(1101000101)$ and $w=(10110)$.  In this case we have $v\circ w = (11\underline{1}1\underline{0}\underline{1}\underline{1}1\underline{0}1)$, where the underlined terms indicate where we have inserted $w$ into $v$. Then $c(v,w) = 8$. We interpret $c(v,w)$ as the number of crossings in a certain diagram associated to $(v,w)$.  First, draw the shuffle permutation associated to $v$:
\[
\begin{minipage}{1.6in}
\vskip.1in
\labellist
\small
\pinlabel $1$ at 12 60
\pinlabel $1$ at 24 60
\pinlabel $0$ at 36 60
\pinlabel $1$ at 48 60
\pinlabel $0$ at 60 60
\pinlabel $0$ at 72 60
\pinlabel $0$ at 84 60
\pinlabel $1$ at 96 60
\pinlabel $0$ at 108 60
\pinlabel $1$ at 120 60
\endlabellist
\begin{center}\includegraphics[scale=.8]{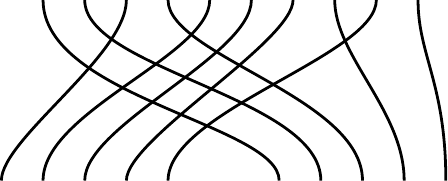}\end{center} 
\end{minipage}
\]
The 1's in the sequence $w$ tell us which strands corresponding to zeroes of $v$ are ``on.''  We will color green and red the strands which are ``on'' and ``off,'' respectively, in the case of $w=(10110)$:
\[
\begin{minipage}{1.6in}
\vskip.1in
\labellist
\small
\pinlabel $1$ at 12 60
\pinlabel $1$ at 24 60
\pinlabel $\textcolor{mygreen}{0}$ at 36 60
\pinlabel $1$ at 48 60
\pinlabel $\textcolor{myred}{0}$ at 60 60
\pinlabel $\textcolor{mygreen}{0}$ at 72 60
\pinlabel $\textcolor{mygreen}{0}$ at 84 60
\pinlabel $1$ at 96 60
\pinlabel $\textcolor{myred}{0}$ at 108 60
\pinlabel $1$ at 120 60
\endlabellist
\begin{center}\includegraphics[scale=.8]{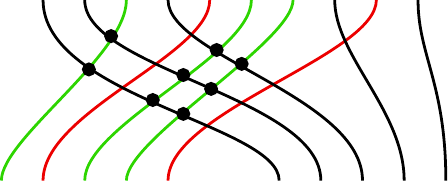}\end{center} 
\end{minipage}
\]
Then $c(v,w)$ is the total number of crossings between the 1-strands and the ``on'' 0-strands, shown here as dots.
\end{example}
	
We now prove:

\begin{lemma}\label{lemma:HH0recursion}
Let $v\in \{0,1\}^n$ be a sequence with $\ell_0$ ones and $k_0=n-\ell_0$ zeroes.  The functions $f_v(q,0,t)$ satisfy the recursion $f_{\emptyset}(q,0,t)=1$ and
\[
f_v(q,0,t) = \sum_{w\in \{0,1\}^{k_0}} t^{\binom{\ell_0}{2}+c(v,w)} q^{k_1} f_w(q,0,t) 
\]
where $w$ has $k_1$ zeroes.
\end{lemma}
\begin{proof}
Setting $a=0$ in Definition \ref{defn-gradings}, we see that $\PC_{v,w}(q,0,t)$ is the product of $t^{\ell(i)+m(i)}$ over all indices $i$ such that $v_i=1$.  The $t^{\ell(i)}$ factors contribute $t^{\binom{\ell_1}{2}}$ and the $t^{m(i)}$ factors contribute $t^{c(v,w)}$.  Thus, setting $a=0$ in the alternate recursion Proposition \ref{prop:altRecursion} gives the statement.  We are also using the result of Proposition \ref{prop-rule1} in the special case where $v$ is the zero sequence.
\end{proof}

\begin{proof}[Proof of Theorem \ref{thm-degreezerosolve-intro}] Consider a sequence $\vb=(v^{(0)},v^{(1)},\ldots,)$ of sequences $v^{(r)}\in \{0,1\}^n$ such that $v^{(r)}_i\leq v^{(r+1)}_i$ for all $i\in\{1,\ldots,n\}$ and all $r\geq 0$.  Assume that $v^{(r)}=(1,1,\ldots,1)$ for $r\gg 0$.   Then $\vb$ defines a function $\sigma:\{1,\ldots,n\}\rightarrow \Z_{\geq 0}$, where $\sigma(i)$ is the smallest $r$ such that $v^{(r)}_i=1$.  It is easy to see that this yields a bijection between sequences $\vb$ and functions $\sigma$.

Let $k(r)$ denote the number of zeroes of $v^{(r)}$, and let $w^{(r)}\in \{0,1\}^{k(r)}$ denote the sequence such that $w^{(0)}=v^{(0)}$ and $v^{(r-1)}\circ w^{(r)}=v^{(r)}$ for $r\geq 1$.  Clearly this establishes a bijection between increasing sequences $\vb$ and eventually empty sequences $\wb=(w^{(1)},w^{(2)},\ldots )$ such that $(w^{(r)},w^{(r+1)})$ are compatible for all $r$.
 
With these notions in place, we apply Lemma \ref{lemma:HH0recursion} to $f_{(00\cdots 0)}$ iteratively. After one application, we see that $f_{00\cdots 0}$ is a sum over sequences $w^{(0)}$ of $f_{w^{(0)}}$, weighted by some monomials in $q$ and $t$.  We apply the recursion again.  The result can be viewed as a sum over sequences $w^{(0)}$ and $w^{(1)}$ with $(w^{(0)},w^{(1)})$ compatible.  Iterating indefinitely, we see that $f_{00\cdots 0}$ can be expressed as a sum over all sequences $\wb$ of some monomials in $q$ and $t$.  Using the bijection between the $\wb$'s and the $\vb$'s, we regard this as a sum over all sequences $\vb$ of some  monomials $g_{\vb}$ computed from $\vb$.  We claim that $g_{\vb}=t^{a(\sigma)+b(\sigma)}q^{c(\sigma)}$, which would prove the theorem.

Suppose an index $i$ is such that $v^{(r-1)}=0$ but $v^{(r)}=1$, with $r\geq 1$.  Then the $t$ contribution at the $r$-th step is $t^{c(w^{(r-1)},w^{(r)})+\binom{\sigma\inv(r)}{2}}$ by Lemma \ref{lemma:HH0recursion}.  Taking the product over all $r\geq 1$ accounts for the factor of $t^{a(\sigma)+b(\sigma)}$.

Finally, the number of zeroes in $w^{(r)}$ of indices $i$ such that $\sigma(i)>r$, that is, $\sum_{s>r}|\sigma\inv(s)|$.  Each of these contributes a factor of $q$.  Taking the product over all $r\geq 0$ yields $q$ to the power of
\[
\sum_{0\leq r<s}|\sigma\inv(s)| = \sum_{0\leq s}s|\sigma\inv(s)|=\sum_{i=1}^n \sigma(i)=c(\sigma).
\]
This accounts for the factor of $q^{c(\sigma)}$. \end{proof}

%% file: FTHHHcomputations.tex
\appendix

\section{Miscellaneous computations}
\label{sec:appendix}
In this appendix we illustrate the usefulness of our method with a few computations of triply graded homology for certain torus knots.  For the reader's convenience we present our results with the proper normalization, and we state how to obtain classical invariants from them.

Let $\PC_\b(Q,A,T)$ denote the Poincar\'e series of $\HHH(F(\b))$, where $Q,A,T$ denote the usual quantum degree, homological degree, and Hochschild degree, respectively.   This $\PC_\b$ is an invariant of the braid closure $L=\hat{b}$ up to multiplication by a unit in $\Z[A^{\pm},Q^\pm,T^\pm]$.  The precise normalization which yields a link invariant requires that we introduce half-integral powers of $A$ and $T$:
\[
\PC_L(Q,A,T)=T^{-e(\b)}Q^n (Q^{-1}A^{1/2}T^{1/2})^{e(\b)-n} \PC_\b(Q,A,T),
\]
where $e(\b)$ is the braid exponent (signed number of crossings).  The decategorification corresponds to specializing $T=-1$.  To avoid choosing a square root of $-1$, we first rewrite $\PC_L$ in terms of $Q,T$, and the Homfly variable $\a=A^{1/2}T^{1/2}Q\inv$:
\[
\PC_L(Q,\a,T):=T^{-e(\b)}Q^n \a^{e(\b)-n} \PC_\b(Q,\a,T).
\]
We call $\PC_L$ the \emph{super polynomial}.  Setting $T=-1$ recovers the Homfly polynomial in variables $\a, Q$.  The $\sl_N$ specialization is then obtained by setting $\a=Q^N$.  For reference, the invariant of the unknot is
\[
\PC_U(Q,\a,T) = \frac{\a\inv + \a T\inv}{Q\inv - Q}.
\]

Recall that we prefer the variables $t=T^2Q^{-2}$, $q=Q^2$, $a=A Q^{-2}$.  Thus, we will usually rewrite $P_L$ in terms of these variables.  The decategorification is obtained by setting $t^{1/2}=-q^{-1/2}$; if one wishes to avoid working in a ring with $\sqrt{-1}$, then one should also set $a^{1/2}(tq)^{1/4}=-\a$.  For knots, it turns out that the \emph{reduced superpolynomial} $\tilde P_L(q,a,t):= P_L(q,a,t)/P_U(q,a,t)$ is a Laurent polynomial in $q^{1/2},a,t^{1/2}$, so no technical issue arises from the decategorification $t^{1/2}\mapsto -q^{-1/2}$.  The $\sl_N$ specialization is then obtained by setting $a=-q^N$.  The following computations were all done by hand.  We omit their derivations, in the interest of readability and length.

\begin{example}\label{ex:2nTorusKnots}
The reduced superpolynomial of the $(2,2k+1)$ torus knot is
\[
a^k (tq)^{-k/2}\bigg(t^k+qt^{k-1}+\cdots + q^k + a(t^{k-1}+qt^{k-2}+\cdots + q^{k-1})\bigg)
\]
In particular, the superpolynomial of the right-handed trefoil---that is, the $(2,3)$ torus knot---is $a(tq)^{-1/2}(q+t+a)$.  The decategorification is $-a(q+q\inv+a)$, and the $\sl_N$ specialization is $q^{N-1}+q^{N+1}-q^{2N}$.
\end{example}

\begin{example}
The reduced superpolynomial of the $(3,4)$ torus knot is
\[
a^3(tq)^{-3/2}\bigg(t^3+qt^2+qt+q^2t+q^3+a(t^2+t+qt+q+q^2)+a^2\bigg).
\]
The decategorification is
\[
-a^3\Big(q^{-3}+q^{-1}+1+q+q^3+a(q^{-2}+q^{-1}+1+q+q^2)+a^2\Big).
\]
Setting $a=-q^2$ and $tq=1$ gives the correct $\sl_2$ specialization (Jones polynomial):
\[
q^{6}\Big(q^{-3}+q^{-1}+1+q+q^3- (1+q+q^2+q^3+q^4)+q^4)\Big) = q^3+q^5-q^8.
\]
\end{example}

\begin{example}
The reduced superpolynomial of the $(3,5)$ torus knot is $a^4 (tq)^{-2}$ times 
\begin{eqnarray*}
&& t^4+qt^3+qt^2+q^2t^2+q^2t+q^3t+q^4 + a\Big(t^3+t^2+qt^2+2qt+q^2t+q^2+q^3\Big) +a^2\Big(q+t\Big).
\end{eqnarray*}
\end{example}

\begin{example}
The reduced superpolynomial of the $(4,5)$ torus knot is
\begin{eqnarray*}
&& a^6 (tq)^{-3} \bigg(t^6 + q t^5 + q t^4 + q t^3 + q^2 t^4 + q^3 t^2 + q^2 t^3 +
   q^2 t^2 + q^3 t + q^3 t^3 + q^4 t^2 + q^4 t + q^5 t + q^6\\
&& \ \ \ \ \ +\ \  a\Big(t^5 + t^4 + t^3 + q t^4 + q^2 t^3+ 2 q t^3 + 2q t^2 + q t + 2 q^2 t^2  +2 q^2t+ 2 q^3 t + q^3 t^2  \\
&& \ \ \ \ \ +\ \     q^4 t + q^3+ q^4 + q^5\Big) + a^2\Big(t^3 + t^2 + t + qt^2 + qt + q^2 t + q + q^2 + q^3\Big) + a^3 \bigg).
\end{eqnarray*}
\end{example}

\begin{observation}
Each of the above polynomials is symmetric with respect to exchanging $q$ and $t$.  Further, the smallest $a$-degree summands of the Poincar\'e series of the $(n,n+1)$ torus knots are the $q,t$ Catalan numbers, for $n=2,3,4$.  This verifies a conjecture in \cite{GORS} in these cases.
\end{observation}